\documentclass[preprint,12pt]{elsarticle}
\usepackage{amsmath,amssymb,amsthm}
\usepackage{booktabs}
\usepackage{subeqnarray} 
\usepackage{enumitem}
\usepackage{graphicx}
\usepackage{algorithm}
\usepackage{algpseudocode}
\usepackage[colorlinks=true,linkcolor=blue,citecolor=blue]{hyperref}

\journal{Computer Methods in Applied Mechanics and Engineering}

\newtheorem{theorem}{Theorem}[section]
\newtheorem{lemma}[theorem]{Lemma}
\newtheorem{proposition}[theorem]{Proposition}
\newtheorem{corollary}[theorem]{Corollary}
\theoremstyle{definition}
\newtheorem{definition}[theorem]{Definition}
\newtheorem{assumption}[theorem]{Assumption}
\theoremstyle{remark}
\newtheorem{remark}[theorem]{Remark}
\newcommand{\Reals}[1]{{\rm I\! R}^{#1}}

\newcommand{\jump}[1]{[\![#1]\!]}
\newcommand{\avg}[1]{\{\!\!\{#1\}\!\!\}}
\newcommand{\Th}{\mathcal{T}_h}
\newcommand{\Eh}{\mathcal{E}_h}
\newcommand{\EI}{\mathcal{E}_h^{I}}

\newcommand{\Wn}{W_n}
\newcommand{\Enr}{E_h}
\newcommand{\V}{V_{n,h}}
\newcommand{\Q}{Q_h}
\newcommand{\Qo}{\mathring Q_h}
\newcommand{\bseps}{{\epsilon}}

\newcommand{\enorm}[1]{\|#1\|_{\mathcal E}}
\newcommand{\divh}{\operatorname{div}_h}

\begin{document}

\begin{frontmatter}

\title{A locking-free enriched linearized neural network for linear
elasticity}

\author[ku]{Bonyong Koo}
\author[txst]{Young Ju Lee\corref{cor1}}
\ead{yjlee@txstate.edu}
\cortext[cor1]{Corresponding author.}
\affiliation[ku]{organization={Department of Mechanical Engineering, Kunsan National University},
            city={Gunsan},
            country={Republic of Korea}}
\affiliation[txst]{organization={Department of Mathematics, Texas State University},
            city={San Marcos},
            state={TX},
            country={USA}}

\begin{abstract}Volumetric locking remains a fundamental challenge for
neural-network Galerkin methods in nearly incompressible elasticity.
We introduce the \emph{locking-free enriched linearized neural network}
(LF-ENN), to the best of our knowledge, the first neural-network
Galerkin method with a rigorous locking-free guarantee.
Building on the enriched Galerkin framework, the method combines
neural approximation with a structure-preserving enrichment that
reduces the discrete incompressibility constraint to a graph-Laplacian
problem. This construction enables divergence-free corrections
and error estimates uniform in the Lam\'e parameter
$\lambda\in(0,\infty]$.
Nitsche boundary enforcement and tailored stabilization overcome
the lack of conventional finite element trace inequalities for
neural trial spaces.
Under suitable approximation and mesh-resolution assumptions,
we establish stability and error estimates independent of network
width, mesh size, and $\lambda$, including an $O(h^s)$ estimate
for reduced-regularity solutions with $s>1/2$.
Numerical experiments in two and three dimensions confirm
locking-free behavior and demonstrate substantial gains in
displacement accuracy.
In particular, LF-ENN achieves approximately twelve times smaller
displacement error than a higher-order locking-free enriched
Galerkin method while using fewer than one-third as many unknowns.
These results establish a rigorous neural Galerkin framework
that combines incompressibility robustness with markedly
improved approximation accuracy.
\end{abstract}

\begin{highlights}
\item A locking-free neural-network Galerkin method for nearly incompressible linear elasticity. 
\item Elementwise dilation enrichment yields a discrete divergence represented by a dual-graph Laplacian. 
\item A displacement-only formulation for finite Lam\'{e}  parameters extends to the incompressible limit through a mixed formulation.
\item Uniform stability and error estimates are established without network-dependent inverse trace inequalities.
\item Numerical experiments demonstrate locking-free accuracy in two and three dimensions, including Cook's membrane.
\end{highlights}


\end{frontmatter}

\section{Introduction}\label{sec:intro}
For nearly incompressible materials the Lam\'e parameter $\lambda$ is much larger than the shear modulus $\mu$.
The divergence then carries the pressure: $p=\lambda\nabla\cdot u$ stays bounded as $\lambda\to\infty$, so
$\nabla\cdot u=O(\lambda^{-1})$, and a displacement method must keep $\lambda\nabla\cdot U$ close to $p$. This requires
enough discrete displacements with zero discrete divergence to approximate divergence-free displacements;
otherwise the term $\lambda\|\nabla\cdot U\|_0^2$ pulls $U$ away from $u$, and the error grows with $\lambda$
(volumetric locking). For continuous $P^1$ with one constraint per element, the constraints are about as many as
the unknowns, so almost no field is discretely divergence-free. Linearized neural networks, whose inner parameters
are fixed so that the discrete problem is linear \cite{LMX25}, approximate smooth functions at high order with few
unknowns, but they have no functions supported on a single element. Consequently, a plain network Galerkin
method using the volumetric term
$\lambda(\nabla\cdot U,\nabla\cdot v)$
may suffer from volumetric locking.
As illustrated in Table~1, its displacement
error increases as $\lambda$ grows,
whereas the LF-ENN displacement error
remains nearly unchanged.

Due to its significance, there have been a number of studies on this locking issue and developments of locking-free methods. Locking of conforming elements was analysed by Babu\v{s}ka and Suri \cite{BaSu92}. Classical remedies
are nonconforming elements \cite{Fal91,BS92,Bre04}, inf-sup stable mixed methods \cite{BF91}, reduced integration, which is
equivalent to a mixed method with discontinuous pressure \cite{MH78} (as is our method, Section~\ref{sec:mixedform}), and
discontinuous Galerkin methods \cite{HL02,CSW06}. The enrichment of \cite{YLZ22} also yields an EG method for Stokes
\cite{YHLA22}, made pressure-robust in \cite{HLMY22}. Penalty-free variants include an arbitrary-order extension of \cite{YLZ22}
in a modified weak Galerkin framework \cite{STZZ24} and face enrichments with a piecewise constant weak
divergence \cite{SYZ25}. On the other hand, to the best of our knowledge, this is the first neural-network Galerkin discretization for linear elasticity equipped with a rigorous locking-free error analysis uniform in the Lam\'e parameter. There are many neural methods to solve solid mechanics problem. Neural methods for solid mechanics include PINNs and the deep energy method \cite{HRMGJ21,SAGN20}; for nearly
incompressible materials PINNs show a locking-like instability \cite{DKLMP25}. Linearized networks include extreme learning
machines \cite{DL21}, random features \cite{CCEY22}, finite neuron and ReLU$^k$ spaces \cite{Xu20,LMX25} and Galerkin neural
networks \cite{AD21}. 

The main purpose of this paper is to develop a neural network based Galerkin method which is locking-free. Our idea is based on the
enriched Galerkin (EG) finite element methods
\cite{SL09,LLW16}. Among others, Yi, Lee and Zikatanov \cite{YLZ22} showed that one dilation field $c_T(x-x_T)$ per element enrichment makes
continuous $P^1$ locking-free (LF-EG), with a penalty term $\lambda^2\beta_{\mathrm{EG}}\sum_eh_e\langle\jump{\nabla\cdot U},\jump{\nabla\cdot v}\rangle_e$
on the jumps of the divergence. Our strategy is to replace
the conforming part by a linearized neural network $\Wn$ and keep the dilation fields $\Enr$. Thus, this method can be considered to be hybrid, similar to the enriched linearized neural network (ENN) method \cite{KJL26}, which is locally conservative. On the other hand, in the discrete bilinear form, 
we employ the piecewise constant averaged-trace divergence instead of the element-wise 
divergence operator. The averaged-trace discrete divergence operator allows us to eliminate the divergence-jump penalty term employed in
\cite{YLZ22}, while retaining the stability properties required for locking-free error analysis. Moreover, the
underlying construction is not restricted to piecewise linear finite element functions or neural-network spaces;
it can be extended to other conforming approximation
spaces, provided that the requisite approximation and stability conditions are satisfied.

The main contributions of this paper can be summarized
as follows:

\begin{enumerate}[wide=0pt,label=\roman*]
\item[I:] \emph{A neural-network Galerkin method that is proved to be locking-free.} 
The displacement and pressure error estimates are independent
of the Lam\'e parameter $\lambda\in(0,\infty]$ on general
shape-regular meshes (Theorem~\ref{thm:main},
Corollary~\ref{cor:anymesh}), including solutions with reduced
regularity (Theorem~\ref{thm:lowreg}).

\item[II:] \emph{One enrichment unknown for each incompressibility constraint.} The discrete divergence on the dilation enrichment is
represented by a dual-graph Laplacian (Lemma~\ref{lem:graph}).
Consequently, every conforming field can be corrected to be
discretely divergence-free (Proposition~\ref{prop:kernel}),
with uniform approximation control under the inf-sup condition.

\item[III:] \emph{No boundary factor or network-dependent inverse trace constant.} Nitsche's method imposes Dirichlet conditions weakly, while
shear least-squares and Hessian stabilizations provide
continuity without a network-dependent inverse trace
inequality (Lemma~\ref{lem:cont}).
Uniform inf-sup stability follows when the network resolves
the mesh (Lemma~\ref{lem:infsup}).

\item[IV:] \emph{High accuracy per unknown.} LF-ENN outperforms the compared Taylor-Hood and enriched
Galerkin methods in the reported benchmarks
(Section~\ref{sec:numerics}).
For Cook's membrane, the tip displacement agrees with
the reference value within $0.02\%$ using 914 unknowns
(Section~\ref{sec:cook}).
\end{enumerate}

Throughout the paper, let $n$ denote the width of the network (its number of neurons) and $M:=\#\Th$ the
number of elements of the auxiliary mesh $\Th$. The symbol $\nu_T$ denotes the outward unit normal of $T \in \Th$. We use the standard notation for Sobolev spaces, such as $H^\ell(\Omega)$ with $H^0(\Omega)$ being $L^2(\Omega)$, and
$H^1_{\Gamma_D}(\Omega):=\{v\in H^1(\Omega):\ v=0\text{ on }\Gamma_D\}$. We also use broken Sobolev spaces. For example, $H^2(\Th)$ denotes
the piecewise $H^2$ functions. 

Finally, we use $C$ and $c$ to denote
generic positive constants independent of
$h$, $n$, and the Lam\'e parameter $\lambda$.
These constants may depend on $d$, the
fixed domain $\Omega$, the boundary
partition $(\Gamma_D,\Gamma_N)$, the mesh
shape-regularity constant, and, where
indicated, the quasi-uniformity constant
and stabilization parameters.
We write $a\lesssim b$ when $a\le Cb$.

The rest of the paper is organized as follows. \S~\ref{sec:model} states the model problem and defines the method, \S~\ref{sec:analysis}
proves stability and error estimates, \S~\ref{sec:numerics} and \S~\ref{sec:cook} report numerical experiments, and \S~\ref{sec:conclusions}
discusses open problems.

\section{Linear Elasticity and LF-ENN}\label{sec:model}
Let $\Omega\subset\Reals{d}$ be a bounded, connected Lipschitz polyhedral domain with $\partial\Omega=\Gamma_D\cup\Gamma_N$, $|\Gamma_D|>0$; the
clamped case $\Gamma_D=\partial\Omega$ is included. Given $f\in L^2(\Omega)^d$ and $g\in L^2(\Gamma_N)^d$, we seek a displacement $u$ and pressure $p$ satisfying 
\begin{equation}\label{eq:model}
\begin{aligned}
-\nabla\cdot\big(2\mu\bseps(u)\big)-\nabla p=f,\quad \nabla\cdot u-\lambda^{-1}p=0\quad&\text{in }\Omega,\\
u=0\quad\text{on }\Gamma_D,\qquad \big(2\mu\bseps(u)+p\mathbf I\big)\nu=g\quad&\text{on }\Gamma_N,
\end{aligned}
\end{equation}
where $ 0 < \mu$, $0<\lambda\le\infty$, and $\bseps(u)$ is the linearized strain tensor.
For $0<\lambda<\infty$, we define the pressure $p:=\lambda\nabla\cdot u$; then \eqref{eq:model} is the linear
elasticity. For $\lambda=\infty$ the second equation of
\eqref{eq:model} reads $\nabla\cdot u=0$, and $p$ is the Lagrange multiplier of this constraint (the Stokes equation), with $p\in L^2_0(\Omega)$
if $\Gamma_N=\emptyset$. For convex polygons and a clamped boundary, $\|u\|_2+\|p\|_1 \lesssim \|f\|_0$ uniformly in $\lambda$ \cite{BS92};
where $\Gamma_D$ meets $\Gamma_N$ the solution is in general less regular, and we give an estimate for $u\in H^{1+s}$, $p\in H^s$,
$\frac12<s<1$ (Theorem~\ref{thm:lowreg}).

\subsection{Auxiliary mesh, dual graph and the LF-ENN method}\label{sec:method}
The method uses an auxiliary mesh, on which the enrichment, the discrete divergence and the pressure are defined. Unlike in
ENN \cite{KJL26}, where the auxiliary mesh serves only local conservation, here it also resolves the pressure and therefore
affects the accuracy (Corollary~\ref{cor:rates}). We first introduce the mesh, its dual graph and the related notation.

\subsubsection{Auxiliary mesh and notation}
Let $\Th$ be a shape-regular partition of $\Omega$ into simplices, with $h_T:=\operatorname{diam}T$, $h:=\max_Th_T$, and faces $\EI$
(interior), $\mathcal E_h^D$ (on $\Gamma_D$) and $\mathcal E_h^N$ (on $\Gamma_N$), $h_e:=\operatorname{diam}e$; we assume that every boundary
face lies in $\overline{\Gamma_D}$ or in $\overline{\Gamma_N}$, and that $\Th$ has more than one element. The faces that carry jumps are
$\Eh:=\EI\cup\mathcal E_h^D$. As in ENN \cite{KJL26}, for an interior face $e=\partial T^+\cap\partial T^-$ with
$\nu^\pm$ the unit outward normal of $T^\pm$ and $v^\pm:=v|_{T^\pm}$, we define the average and the jumps
\begin{equation}
\avg v:=\tfrac12\big(v^++v^-\big),\qquad \jump v:=v^+\otimes \nu^++v^-\otimes \nu^-,\qquad \jump v_n:=v^+\cdot \nu^++v^-\cdot \nu^-
\end{equation}
for vector-valued $v$, and $\jump q:=q^+\nu^++q^-\nu^-$ for scalar $q$; on a boundary face, $\avg v:=v$, $\jump v:=v\otimes\nu$,
$\jump v_n:=v\cdot\nu$ and $\jump q:=q\nu$, with $\nu$ the outward normal of $\Omega$. Since $\nu^-=-\nu^+$, with $\nu_e:=\nu^+$,
\[
\avg\sigma:\jump v=\avg{\sigma}\nu_e\cdot(v^+-v^-),\qquad |\jump v|=|v^+-v^-|,\qquad \jump v_n=(v^+-v^-)\cdot \nu_e
\]
for symmetric tensors $\sigma$; on faces where only a normal is needed, $\nu_e$ denotes $\nu^+$ or, on $\partial\Omega$, the outward normal.
Broken operators act element by element: $\bseps_h(v)|_T:=\bseps(v|_T)$, $(\nabla_h\!\cdot v)|_T:=\nabla\cdot(v|_T)$, and
$(\cdot,\cdot)$ is the $L^2(\Omega)$ product. For face functions, $\langle\phi,\psi\rangle_{\Eh}:=\sum_{e\in\Eh}\int_e\phi\cdot\psi\,ds$, with
$\phi:\psi$ in place of $\phi\cdot\psi$ for matrices. For $k\ge0$, $Q_k$ denotes the $L^2(\Omega)$ projection onto the discontinuous
piecewise polynomials of degree $k$ on $\Th$. We use four elementary facts repeatedly.
\begin{enumerate}[wide=0pt]
\item[(F1)] If $T$ is a simplex and $e$ a face of $T$ with height $H_e$ (the distance from the opposite vertex to the
hyperplane of $e$), then $(x-x_T)\cdot \nu_T=H_e/(d+1)$ for all $x\in e$, and $|T|=|e|H_e/d$. Hence
$\int_e(x-x_T)\cdot \nu_T\,ds=\frac d{d+1}|T|$.
\item[(F2)] For $T\in\Th$, a face $e\subset\partial T$ and $\phi\in H^1(T)$,
$\|\phi\|_{0,e}^2\lesssim h_T^{-1}\|\phi\|_{0,T}^2+h_T|\phi|_{1,T}^2$.
\item[(F3)] Let $\chi_D:=\sum_x\varphi_x$ be the sum of the continuous piecewise linear hat functions of the vertices $x$
of $\Th$ on $\overline{\Gamma_D}$, $\Omega_D$ the union of the elements that have such a vertex, and $\mathcal E_h^{J}\subset\mathcal E_h^N$ the
traction faces with a vertex on $\overline{\Gamma_D}$. Then $\chi_D=1$ on $\Gamma_D$, $0\le\chi_D\le1$, $\operatorname{supp}\chi_D\subset\overline{\Omega_D}$,
$|\nabla\chi_D|\lesssim h_T^{-1}$ on $T$, $\chi_D=0$ on the traction faces not in $\mathcal E_h^J$, and for $\phi\in H^1(\Omega)$,
\[
\sum_{T\subset\Omega_D}h_T^{-2}\|\phi\|_{0,T}^2+\sum_{e\in\mathcal E_h^J}h_e^{-1}\|\phi\|_{0,e}^2\lesssim\sum_{e\in\mathcal E_h^D}h_e^{-1}\|\phi\|_{0,e}^2+\|\nabla\phi\|_{0,\Omega_D}^2 .
\]
\item[(F4)] If $\sigma\in H(\mathrm{div};\Omega)^d$ is piecewise $H^1$ and $\zeta$ is piecewise $H^1$, then
\[
\sum_{e\in\EI}\int_e\sigma\nu_e\cdot(\zeta^+-\zeta^-)\,ds+\int_{\partial\Omega}\sigma\nu\cdot\zeta\,ds=(\sigma,\nabla_h\zeta)+(\nabla\cdot\sigma,\zeta).
\]
\end{enumerate}
(F1) holds because $x_T$ is the barycentre, which lies at distance $H_e/(d+1)$ from the hyperplane of every face. For each Dirichlet vertex $x$, let $\omega_x$ denote its
vertex patch, define
\[
h_x:=\max_{T\subset\omega_x}h_T,
\]
and choose a Dirichlet face $e_x\in\mathcal E_h^D$
containing $x$. By shape regularity and conformity,
the vertex patch consists of a uniformly bounded number
of face-connected elements with comparable diameters;
in particular, $h_{e_x}\simeq h_x$.
The scaled Friedrichs inequality on $\omega_x$ yields
\[
h_x^{-2}\|\varphi\|_{0,\omega_x}^2
\lesssim
h_{e_x}^{-1}\|\varphi\|_{0,e_x}^2
+\|\nabla\varphi\|_{0,\omega_x}^2.
\]
Summing over the Dirichlet vertices and using bounded
overlap proves the element contribution in (F3).
The face contribution follows from (F2). (F4) is elementwise integration by parts,
$\sum_T\int_{\partial T}\sigma\nu_T\cdot\zeta|_T=\sum_T[(\sigma,\nabla\zeta)_T+(\nabla\cdot\sigma,\zeta)_T]$, together with the continuity of $\sigma\nu_e$ across
interior faces. Finally, since $|\Gamma_D|>0$, Korn's inequality $\|\nabla\phi\|_0\lesssim\|\bseps(\phi)\|_0+\|\phi\|_{0,\Gamma_D}$ holds for
$\phi\in H^1(\Omega)^d$, and the broken Korn inequality $\|\nabla_h v\|_0\lesssim\|\bseps_h(v)\|_0+(\sum_{e\in\Eh}h_e^{-1}\|\jump v\|_{0,e}^2)^{1/2}$ holds for
piecewise $H^1$ fields \cite{Bre04,HLX25}.

\subsubsection{Dual graph of $\mathcal{T}_h$} 
We shall index the elements as $\Th=\{T_1,\dots,T_M\}$, $M:=\#\Th$. The dual graph of $\Th$ has one node per element, which may be placed
at its barycentre $x_{T_i}$, and one face per interior face, joining the two elements that share it; we write $T_j\sim T_i$ if $T_i$
and $T_j$ share an interior face. The affinity (adjacency) matrix $\mathcal W\in\Reals{M\times M}$ and the degree matrix $\mathcal D$ of
the dual graph are
\begin{equation}
\mathcal W_{ij}:=\begin{cases}1,&T_j\sim T_i,\\0,&\text{otherwise},\end{cases}\qquad
\mathcal D:=\operatorname{diag}\Big(\sum_{j=1}^M\mathcal W_{ij}\Big)_{i=1}^M,
\end{equation}
so that $\mathcal D_{ii}$ is the number of interior faces of $T_i$. The graph Laplacian of the dual graph is
\begin{equation}
L:=\mathcal D-\mathcal W,\qquad\text{i.e.}\quad(Lz)_i=\sum_{j:\,T_j\sim T_i}(z_i-z_j),\quad z\in\Reals{M}.
\end{equation}
It is symmetric positive semidefinite, $z^\top Lz=\sum_{T_i\sim T_j,\,i<j}(z_i-z_j)^2$ with one term per interior face, and for
connected $\Omega$ its kernel consists of the constant vectors. For a vector indexed by elements we write $z_T:=z_i$ and
$(Lz)_T:=(Lz)_i$ when $T=T_i$.

\subsubsection{Trial space, forms and LF-ENN}
The trial space $\V$ consists of a conforming part $\Wn$ and an enrichment. The conforming part is a finite-dimensional space with
\begin{equation}\label{eq:W}
\Wn\subset H^2(\Omega)^d .
\end{equation}
No boundary condition is built into $\Wn$. We set
\begin{gather*}
\V := \Wn+\Enr,\qquad \Enr:=\{v:\ v|_T=c_T(x-x_T),\ c_T\in\Reals{}\},
\end{gather*}
where $x_T$ is the barycentre, and the pressure spaces
\begin{gather*}
\Q:=\{q:\ q|_T\in P^0(T)\},\qquad \Qo:=\begin{cases}\Q\cap L^2_0(\Omega),&\Gamma_N=\emptyset,\\ \Q,&\Gamma_N\neq\emptyset.\end{cases}
\end{gather*}
We write $v = v_W + v_E$. Since $v_W$ is continuous, $\jump v=\jump{v_E}$ on interior faces, and $\jump v=(v_W+v_E)\otimes\nu$ on clamped faces.
Our main instance is the linearized neural network: 
\begin{gather*}
\Wn:=\operatorname{span}\{\phi_j\,\mathbf e_i:\ 1\le j\le n,\ 1\le i\le d\},\quad
\phi_j(x) = \sigma_k\big(\theta_j\cdot(x,1)\big),\quad k\ge2,
\end{gather*}
with $\theta_j$ quasi-uniform on $\mathbb S^d$ \cite{LMX25,KJL26}, and the ReLU$^k$ activation $\sigma_k(t):=\max\{0,t\}^k$.
An affine change of the input does not change this class: for an invertible $A\in\Reals{d\times d}$ and $c\in\Reals d$,
$\theta\cdot(Ax+c,1)=(M\theta)\cdot(x,1)$ with $M(\omega,b):=(A^\top\omega,\,c\cdot\omega+b)$, so
$\sigma_k\big(\theta_j\cdot(Ax+c,1)\big)=|M\theta_j|^k\sigma_k\big(\theta_j'\cdot(x,1)\big)$ with $\theta_j':=M\theta_j/|M\theta_j|\in\mathbb S^d$. The map
$\theta\mapsto M\theta/|M\theta|$ is a diffeomorphism of $\mathbb S^d$, so the $\theta_j'$ are again quasi-uniform, with constants depending on $M$. In
particular, the input can be translated and scaled to the domain at hand. For $k\ge2$ the features are $C^1$ and piecewise polynomial, so \eqref{eq:W} holds. The neural resolution is $\delta := n^{-1/d}$. We call the method
LF-ENN for this network.

Let $\alpha>0$ be a penalty parameter. For $u,v\in\V+H^2(\Th)^d$, the shear form is the NIPG form
\begin{equation}\label{eq:Ah}
\begin{aligned}
A_h(u,v) &= 2\mu(\bseps_h(u),\bseps_h(v))-\langle\avg{2\mu\bseps(u)},\jump v\rangle_{\Eh}
+\langle\jump u,\avg{2\mu\bseps(v)}\rangle_{\Eh} \\
&\quad +\mu\alpha\langle h_e^{-1}\jump u,\jump v\rangle_{\Eh}.
\end{aligned}
\end{equation}
The second term is the consistency term, the third makes the form nonsymmetric, and the fourth is the jump penalty. On interior faces
every jump is a jump of the enrichment; on clamped faces $\jump v=v\otimes\nu$, and the boundary terms are the nonsymmetric Nitsche terms
\cite{Nit71} that impose $u=0$ on $\Gamma_D$ weakly, for the network and the dilation fields alike. We use $\alpha=1$.

Networks satisfy no inverse trace inequality, so the face integrals of the network traction in \eqref{eq:Ah} cannot be bounded by
the energy of the network alone. As in ENN \cite{KJL26}, we add a Galerkin least-squares term on the shear part. Let $c_\tau>0$ be fixed,
$\tau_T:=c_\tau h_T^2/\mu$, and
\begin{equation}\label{eq:gls}
G_h(u,v) := \sum_{T\in\Th}\tau_T\big(\nabla\cdot2\mu\bseps(u),\nabla\cdot2\mu\bseps(v)\big)_T,\,\,
\ell_h(v):=-\sum_{T\in\Th}\tau_T\big(f,\nabla\cdot2\mu\bseps(v)\big)_T .
\end{equation}
They come from the least-squares term $\sum_T\tau_T(-\nabla\cdot2\mu\bseps(U)-f,-\nabla\cdot2\mu\bseps(v))_T$ of the shear part. They act only on $\Wn$,
because $\bseps(v_E)$ is constant on every element, and they contain no $\lambda$. For the exact solution the residual of the shear part is
$-\nabla\cdot2\mu\bseps(u)-f=\nabla p$; the discrete pressure is piecewise constant and cannot represent it, so the term is consistent only up
to the pressure gradient (Lemma~\ref{lem:consistency}). We use $c_\tau=0.01$.
The shear least-squares term, together with the volume strain energy, permits control of the network traction contributions on interior and Dirichlet faces through integration by parts, without requiring a network-dependent inverse trace inequality, where the jumps of the dilation fields are in the
energy norm (Lemma~\ref{lem:cont}). On the traction faces the dilation fields have no jump term, and the traction of the network must be
controlled directly. Let $\Omega_N$ be the union of the elements that have a vertex on $\overline{\Gamma_N}$ (empty if $\Gamma_N=\emptyset$), $c_H>0$, and
an asymptotically consistent stabilization: 
\begin{equation}\label{eq:hess}
H_h(u,v):=\mu\,c_H\sum_{T\subset\Omega_N}h_T^2\,(\nabla^2u,\nabla^2v)_T .
\end{equation}
It acts only on $\Wn$, since $\nabla^2v_E=0$, and contains no $\lambda$. It is not consistent; its residual is of the size $h|u|_2$, the same order
as the other consistency terms (Lemma~\ref{lem:consistency}). We use $c_H=0.01$ and write $\mathcal A_h:=A_h+G_h+H_h$.

For the divergence form, note that for smooth $p$ elementwise integration by
parts gives $-(\nabla p,v)=(p,\nabla_h\!\cdot v)-\langle p,\jump v_n\rangle_{\Eh}-\langle p,v\cdot\nu\rangle_{\Gamma_N}$ for all $v\in\V$. Replacing $p$ on
the faces of $\Eh$ by the average of a piecewise constant $q$ gives the standard discontinuous Galerkin divergence form \cite{CKSS02,DPE12}: for $q \in Q_h$, 
\begin{equation}
B_h(v,q) := (\nabla_h\!\cdot v,q)-\langle\avg q,\jump v_n\rangle_{\Eh} =\sum_{e\in\EI}\int_e\avg v\cdot\jump q\,ds+\langle q,v\cdot\nu\rangle_{\Gamma_N}, \label{eq:bface}
\end{equation}
the averaged-trace divergence form; the second expression follows from elementwise integration by parts and
$\jump{qv}_n=\avg q\jump v_n+\avg v\cdot\jump q$ on interior faces. We define $\divh:\V\to\Q$ by
$(\divh v,q) = B_h(v,q)$, i.e.\ $(\divh v)|_T=|T|^{-1}\int_{\partial T\setminus\Gamma_D}\avg v\cdot \nu_T\,ds$. If $\Gamma_N=\emptyset$, the interior face
fluxes cancel in $\sum_T|T|(\divh v)|_T$, so $\divh$ maps $\V$ into $\Qo$ in both cases. Thus $\divh v=0$ states local
conservation of volume, i.e., the averaged flux through the boundary of every element vanishes.
\begin{remark}[Conservation of Volume] 
For a displacement $u$, the normal flux $\int_{\partial T}u\cdot\nu_T\,ds=\int_T\nabla\cdot u\,dx$ is
the linearized change of volume of $T$. Accordingly, $|T|\,(\divh v)|_T=\int_{\partial T\setminus\Gamma_D}\avg v\cdot\nu_T\,ds$ is the linearized
change of volume of $T$ when its faces move with the averaged displacement $\avg v$; the average makes this flux single valued, so
that neighbouring elements share each deformed face. The method therefore imposes, element by element,
$\divh U=\lambda^{-1}P$, and for $\lambda=\infty$ local conservation of volume, $\divh U=0$: every element keeps its volume. This differs from the local conservation property in discontinuous Galerkin methods and ENN \cite{KJL26}, which concerns the balance of linear momentum rather than conservation of mass. There, testing the discrete equilibrium equation with the indicator function of an element yields a local balance between the numerical traction flux and the applied body force. In solid mechanics mass is conserved automatically, since every material element keeps its mass and only its density changes, what a displacement method has to control is the change of volume, which is the source of volumetric locking. Here conservation of volume is a constraint on the displacement itself, imposed through $\divh$, and not a balance law obtained from the test functions, thus the momentum balance is not satisfied element by element.
\end{remark} 

The form $B_h$ is the standard discontinuous Galerkin velocity-pressure coupling \cite{CKSS02,HL02,DPE12}, and $\divh$ is the corresponding
discrete divergence \cite{DPE12} with values in the piecewise constants. What is new here is its use as the only term containing $\lambda$,
and its action on the dilation fields (Lemma~\ref{lem:graph}). The load is the standard one,
\begin{equation}\label{eq:loadN}
F_h(v):=(f,v)+\langle g,v\rangle_{\Gamma_N}+\ell_h(v).
\end{equation}

\begin{definition}[LF-ENN, displacement formulation]\label{def:method}
For $0<\lambda<\infty$, find $U\in\V$ such that
\begin{equation}\label{eq:primal}
\mathcal A_h(U,v) + \lambda(\divh U,\divh v) = F_h(v)\qquad\forall v\in\V,
\end{equation} 
where $\divh v$ is the piecewise constant averaged-trace divergence,
\begin{equation}
(\divh v)|_T:=\frac1{|T|}\int_{\partial T\setminus\Gamma_D}\avg v\cdot \nu_T\,ds, \quad\forall T\in\Th. 
\end{equation} 
The method has displacement unknowns only; the pressure is computed after solving,
element by element, as $P=\lambda\divh U$.

\end{definition}

\subsubsection{Mixed formulation}\label{sec:mixedform}
Since $\divh U\in\Q$, $P=\lambda\divh U$ is equivalent to $B_h(U,q)-\lambda^{-1}(P,q)=0$ for all $q\in\Q$, and $(\divh v,P)=B_h(v,P)$. Hence $U$ solves
\eqref{eq:primal} if and only if $(U,P)\in\V\times\Qo$ solves
\begin{equation}\label{eq:method}
\mathfrak B\big((U,P),(v,q)\big)=F_h(v)\qquad\forall(v,q)\in\V\times\Qo,
\end{equation}
where
\[
\mathfrak B\big((w,r),(v,q)\big):=\mathcal A_h(w,v)+B_h(v,r)-B_h(w,q)+\lambda^{-1}(r,q).
\]
For $\lambda=\infty$, \eqref{eq:method} with $\lambda^{-1}:=0$ is the definition. The mixed form is used in the analysis and for $\lambda=\infty$.

\begin{remark}[Other choices]\label{rem:choices}
(a) With $\lambda(\nabla_h\!\cdot U,\nabla_h\!\cdot v)$ the method locks, since divergences of $\Wn$ are not piecewise constant. Replacing the full divergence penalty by the reduced
integration term
\[
\lambda(Q_0\nabla_h\cdot U,Q_0\nabla_h\cdot v)
\]
can alleviate volumetric locking. However, the resulting
formulation has a consistency residual containing
\[
\langle p, [[v]]_n\rangle_{\mathcal E_h},
\]
which is not generally small. The face correction in
$B_h$ removes this leading inconsistency. (b) If $W_n$ is replaced by a conforming finite element
space of fixed polynomial degree on $\mathcal T_h$,
the standard discrete inverse trace inequality is
available, so the shear least-squares stabilization
$G_h$ is not required for continuity.
The corresponding stability and error analysis can
be developed using piecewise $H^2$ regularity and
the appropriate finite element approximation
properties. (c) Let $b\in W^{2,\infty}(\Omega)$ be a fixed
boundary factor satisfying $b=0$ on $\Gamma_D$.
Replacing $W_n$ by $bW_n$ makes the conforming network component vanish on $\Gamma_D$, but does not by itself impose the Dirichlet condition strongly on the full space $bW_n+E_h$: the enrichment fields generally have nonzero boundary traces. Strong imposition would also require constraining or modifying $E_h$. Let $R_hv$ be a $\mathrm{BDM}_1$
reconstruction whose normal face moments agree
with the averaged normal face moments of $v$,
with zero normal flux on $\Gamma_D$. Then
\[
\operatorname{div}R_hv=\operatorname{div}_h v.
\]
Consequently, replacing the load $(f,v)$ by
$(f,R_hv)$ provides a route to pressure robustness:
discretely divergence-free test functions are
reconstructed as exactly divergence-free fields,
so gradient forces are orthogonal to their
reconstructions, subject to the appropriate
boundary conditions and consistency treatment.
The price is that the network must approximate
$u/b$, which may have lower regularity than $u$,
and that continuity of $\mathcal A_h$ may require
network-dependent trace estimates on interior
faces.
\end{remark}

\section{Stability and Error analysis}\label{sec:analysis}
In this section, we shall establish the error analysis. We shall first introduce the energy norm. For $v \in \V + H^2(\Th)^d$,  
\begin{equation}
\enorm v^2:=\|\bseps_h(v)\|_0^2+\alpha|v|_J^2+c_\tau\sum_{T\in\Th}h_T^2\|\nabla\cdot2\bseps(v)\|_{0,T}^2+c_H\sum_{T\subset\Omega_N}h_T^2|v|_{2,T}^2,
\end{equation}
where $|v|_J^2:=\sum_{e\in\Eh}h_e^{-1}\|\jump v\|_{0,e}^2$. For $\eta\in H^2(\Omega)^d$, we define
the mesh-dependent norm
\begin{equation}
\begin{aligned}
\|\eta\|_*
:={}&|\eta|_{1,\Omega}
+\left(
\sum_{T\in\mathcal T_h}
h_T^2|\eta|_{2,T}^2
\right)^{1/2} +\left(
\sum_{T\in\mathcal T_h}
h_T^{-2}\|\eta\|_{0,T}^2
\right)^{1/2}.
\end{aligned}
\label{eq:star-norm}
\end{equation}
Our main result (Theorem~\ref{thm:main}) is the estimate
\begin{equation}
\enorm{u-U}+\mu^{-1}\|Q_0p-P\|_0\lesssim\inf_{w\in\Wn}\|u-w\|_*+h\big(|u|_2+\mu^{-1}|p|_1\big),
\end{equation}
with constants independent of $\lambda$; on quasi-uniform meshes with $n\gtrsim\#\Th$ they are also independent of $n$ and $h$. For solutions
with the reduced regularity of mixed boundary conditions, $u\in H^{1+s}$ and $p\in H^s$ with $\frac12<s<1$, the rate is $h^s$
(Theorem~\ref{thm:lowreg}).
Readers interested mainly in performance may
skip to Section~\ref{sec:numerics}. 

The first result is about the enrichment controlled by the jumps. More precisely, we can establish that 
\begin{lemma}\label{lem:Ejump}
There exists a constant $c_J > 0$, depending
only on the mesh shape-regularity constant,
such that, for every
$v=v_W+v_E\in V_{n,h}$, the following
estimates hold: 
\begin{equation} 
\sum_Tc_T^2|T|\le c_J\sum_{e\in\EI}h_e^{-1}\|\jump v\|_{0,e}^2
\end{equation} 
and 
\begin{align*}
&\|\nabla v_E\|_0 + \Big(\sum_{T \in\mathcal{T}_h}  h_T^{-2}\|v_E\|_{0,T}^2\Big)^{1/2}+\Big(\sum_{e \in \mathcal{E}} h_e^{-1}\|v_E\|_{0,e}^2\Big)^{1/2}\\
&\qquad+\|\nabla v_W\|_0+\Big(\sum_{e\in\mathcal E_h^D}h_e^{-1}\|v_W\|_{0,e}^2\Big)^{1/2}\le C\enorm v. 
\end{align*}
\end{lemma}
\begin{proof}
Let $e=T^+\cap T^-$ be an interior face, $\nu_e=\nu^+$, and $c^\pm=c_{T^\pm}$, $H^\pm$ the
heights of $T^\pm$ over $e$. Then $\jump{v}=\jump{v_E}=g\otimes \nu_e$ with $g:=c^+(x-x_{T^+})-c^-(x-x_{T^-})$, so $\|\jump{v}\|_{0,e}=\|g\|_{0,e}$. By (F1),
$g$ has the constant normal component $g\cdot \nu_e=(c^+H^++c^-H^-)/(d+1)$, and its tangential component is $(c^+-c^-)x_\tau+\text{const}$, an affine
function with slope $c^+-c^-$ in every tangential direction. Since the normal and tangential parts are orthogonal and an affine function
with slope $s$ has $L^2(e)$-norm at least $c\,|s|\,h_e|e|^{1/2}$ after subtracting its mean,
\[
\|g\|_{0,e}^2\ge|e|\Big(\frac{(c^+H^++c^-H^-)^2}{(d+1)^2}+c\,h_e^2(c^+-c^-)^2\Big)\ge c\,|e|h_e^2\big((c^+)^2+(c^-)^2\big).
\]
The last step follows from shape regularity,
which implies $H^\pm/h_e\simeq 1$.
Indeed, the quadratic form
\[
(a,b)\mapsto
\frac{(aH^++bH^-)^2}{h_e^2}
+(a-b)^2
\]
is uniformly positive definite, since its
kernel is trivial and the normalized heights
$H^\pm/h_e$ belong to a compact subset of
$(0,\infty)$. Thus its smallest eigenvalue
is bounded below by a positive constant
depending only on shape regularity. Since $\mathcal T_h$ is a conforming
partition of the connected domain $\Omega$
into more than one element, every element
shares at least one interior face with
another element. Using
$|e|h_e\simeq|T|$ for $e\subset\partial T$, we obtain that  $\sum_Tc_T^2|T|\lesssim\sum_{e\in\EI}h_e^{-1}\|\jump{v}\|_{0,e}^2\le|v|_J^2$. We observe that $\|\nabla v_E\|_0^2=d\sum_Tc_T^2|T|$. Since $|v_E|\le|c_T|h_T$ on $T$, $h_T^{-2}\|v_E\|_{0,T}^2\lesssim c_T^2|T|$ and $h_e^{-1}\|v_E\|_{0,e}^2\lesssim c_T^2|T|$ for $e\subset\partial T$. On a clamped face $\jump v=(v_W+v_E)\otimes\nu$, so $|v_W|\le|\jump v|+|v_E|$ and
$\sum_{e\in\mathcal E_h^D}h_e^{-1}\|v_W\|_{0,e}^2\lesssim|v|_J^2+\sum_Tc_T^2|T|$. Finally $\|\bseps(v_W)\|_0\le\|\bseps_h(v)\|_0+\|\bseps(v_E)\|_0$, and Korn's inequality with
the boundary term on $\Gamma_D$, together with $\|v_W\|_{0,\Gamma_D}^2\le h\sum_{e\in\mathcal E_h^D}h_e^{-1}\|v_W\|_{0,e}^2$, bounds $\|\nabla v_W\|_0$. This completes the proof. 
\end{proof}

We shall now establish the coercivity and continuity. 
\begin{lemma}\label{lem:cont}
For $w,v\in\V$, it holds that $\mathcal A_h(v,v)\ge\mu\enorm v^2$ and
\begin{equation}\label{eq:cont}
|\mathcal A_h(w,v)|\le C_A\,\mu\enorm w\enorm v,
\end{equation} 
where \(C_A\) depends on \(\alpha\), \(c_\tau\), \(c_H\), the mesh shape-regularity constant, and the fixed domain and boundary partition through the relevant Korn and Friedrichs inequalities, but is independent of \(h\), \(n\), and \(\lambda\). For $\eta\in H^2(\Omega)^d$, $v\in\V$ and $q\in\Q$,
\begin{equation}\label{eq:contEta}
|\mathcal A_h(\eta,v)|\le C_A\,\mu\|\eta\|_*\enorm v,\qquad |B_h(\eta,q)|\le C\|\eta\|_*\|q\|_0 .
\end{equation}
\end{lemma}
\begin{proof}
For $w=v$ the second and third terms of $A_h$ cancel, and $\tau_T\mu^2=\mu c_\tau h_T^2$, so we get that 
\begin{eqnarray*}
\mathcal A_h(v,v) &=& 2\mu\|\bseps_h(v)\|_0^2+\mu\alpha|v|_J^2+\mu c_\tau\sum_Th_T^2\|\nabla\cdot2\bseps(v)\|_{0,T}^2 \\
&& +\mu c_H\sum_{T\subset\Omega_N}h_T^2|v|_{2,T}^2\ge\mu\enorm v^2.   
\end{eqnarray*}
This proves the coercivity. Let $w,v\in\V$, $\sigma:=2\mu\bseps(w_W)$, and $I:=\langle\sigma,\jump v\rangle_{\Eh}$. Since $w_W\in H^2(\Omega)^d$,
$\sigma\in H(\mathrm{div};\Omega)^d$, and $\nabla\cdot\sigma=\nabla\cdot2\mu\bseps(w)$ on every element, because $\bseps(w_E)$ is constant there. Thus, we get 
\begin{eqnarray*}
\|\sigma\|_0 + \Big(\sum_Th_T^2\|\nabla\cdot\sigma\|_{0,T}^2\Big)^{1/2} &\le& C\mu\enorm w,\\
\Big(\sum_{e\in\mathcal E_h^N}h_e\|\sigma\nu\|_{0,e}^2\Big)^{1/2} &\le& C\mu\,(1+c_H^{-1/2})\enorm w,
\end{eqnarray*}
the second by (F2) on the elements $T\subset\Omega_N$ that contain a traction face, $h_e\|\sigma\nu\|_{0,e}^2\lesssim\|\sigma\|_{0,T}^2+h_T^2|\sigma|_{1,T}^2$, with
$|\sigma|_{1,T}\le2\mu|w|_{2,T}$ because $\nabla^2w_E=0$, and Lemma~\ref{lem:Ejump}. We shall write $\jump v=\jump{v_E}+\jump{v_W}$, where $\jump{v_W}=0$ on interior faces and
$\jump{v_W}=v_W\otimes\nu$ on clamped faces. By (F4) with $\zeta=v_E$, the sum over $\Eh$ is the full sum minus the traction faces:
\begin{eqnarray*} 
\langle\sigma,\jump{v_E}\rangle_{\Eh} &=& (\sigma,\nabla_hv_E)+(\nabla\cdot\sigma,v_E)-\int_{\Gamma_N}\sigma\nu\cdot v_E\,ds \\ 
&\leq& \|\sigma\|_0\|\nabla v_E\|_0 +  \Big(\sum_Th_T^2\|\nabla\cdot\sigma\|_{0,T}^2\Big)^{1/2}\Big(\sum_Th_T^{-2}\|v_E\|_{0,T}^2\Big)^{1/2}\\
&& + \Big(\sum_{e\in\mathcal E_h^N}h_e\|\sigma\nu\|_{0,e}^2\Big)^{1/2}\Big(\sum_{e\in\mathcal E_h^N}h_e^{-1}\|v_E\|_{0,e}^2\Big)^{1/2}. 
\end{eqnarray*}
Therefore, by Lemma~\ref{lem:Ejump}, it is bounded by $C\mu\enorm w\enorm v$. For the clamped faces let $\zeta:=\chi_Dv_W\in H^1(\Omega)^d$ with $\chi_D$ from (F3). Then $\zeta=v_W$ on $\Gamma_D$,
$\zeta$ vanishes on the traction faces not in $\mathcal E_h^J$, and the divergence theorem on $\Omega$ gives
\[
\langle\sigma,\jump{v_W}\rangle_{\mathcal E_h^D}=\int_{\Gamma_D}\sigma\nu\cdot v_W\,ds=(\sigma,\nabla\zeta)+(\nabla\cdot\sigma,\zeta)-\sum_{e\in\mathcal E_h^J}\int_e\sigma\nu\cdot\zeta\,ds .
\]
Since $\nabla\zeta=\chi_D\nabla v_W+v_W\otimes\nabla\chi_D$, (F3) and Lemma~\ref{lem:Ejump} give
$\|\nabla\zeta\|_0+(\sum_{T\subset\Omega_D}h_T^{-2}\|\zeta\|_{0,T}^2)^{1/2}+(\sum_{e\in\mathcal E_h^J}h_e^{-1}\|\zeta\|_{0,e}^2)^{1/2}\le C\enorm v$, and the three terms are bounded
as before. Hence $|I|\le C\mu\enorm w\enorm v$, without invoking any network-dependent inverse trace inequality.

We shall now consider the continuity on $\V$. For the enrichment part of the traction, $\bseps(w_E)=c_T\mathbf I$ and $h_e\|c_T\nu\|_{0,e}^2\lesssim c_T^2|T|$, so
$|\langle\avg{2\mu\bseps(w_E)},\jump v\rangle_{\Eh}|\le C\mu\enorm w|v|_J$ by Lemma~\ref{lem:Ejump}. The consistency term of $A_h$ is bounded by $C\mu\enorm w\enorm v$. The term $\langle\jump w,\avg{2\mu\bseps(v)}\rangle_{\Eh}$ is bounded in the same way with $w$ and $v$ exchanged, the volume term by $2\mu\enorm w\enorm v$, the penalty by $\mu\alpha|w|_J|v|_J$, and $G_h(w,v)$ and $H_h(w,v)$ by $\mu\enorm w\enorm v$. For the continuity for $\eta\in H^2(\Omega)^d$, we have that $\jump\eta=0$ on interior faces and $\jump\eta=\eta\otimes\nu$ on clamped faces. The volume term is at most
$2\mu|\eta|_1\enorm v$. By (F2), $h_e\|\bseps(\eta)\nu_e\|_{0,e}^2\lesssim|\eta|_{1,T}^2+h_T^2|\eta|_{2,T}^2$ and $h_e^{-1}\|\eta\|_{0,e}^2\lesssim h_T^{-2}\|\eta\|_{0,T}^2+|\eta|_{1,T}^2$,
which bound $\langle\avg{2\mu\bseps(\eta)},\jump v\rangle_{\Eh}$ and the penalty by $C\mu\|\eta\|_*\enorm v$. 

In $\langle\jump\eta,\avg{2\mu\bseps(v)}\rangle_{\Eh}=\int_{\Gamma_D}\eta\cdot2\mu\bseps(v)\nu$, the enrichment part is bounded by
\[
C\mu\Big(\sum_{e\in\mathcal E_h^D}h_e^{-1}\|\eta\|_{0,e}^2\Big)^{1/2}\enorm v,
\]
and the network part, by the divergence theorem with $\zeta=\chi_D\eta$ and $\sigma=2\mu\bseps(v_W)$ exactly as above, by
\[
C\mu\enorm v\Big(\|\nabla(\chi_D\eta)\|_0+\Big(\sum_{T\subset\Omega_D}h_T^{-2}\|\eta\|_{0,T}^2\Big)^{1/2}
+\Big(\sum_{e\in\mathcal E_h^J}h_e^{-1}\|\eta\|_{0,e}^2\Big)^{1/2}\Big)\le C\mu\|\eta\|_*\enorm v,
\]
using (F2) for the last sum. Finally,
\begin{eqnarray*}
|G_h(\eta,v)| &\le& \mu c_\tau^{1/2}\Big(\sum_Th_T^2\|\nabla\cdot2\bseps(\eta)\|_{0,T}^2\Big)^{1/2}\enorm v \lesssim \mu h|\eta|_2\enorm v,\\
|H_h(\eta,v)|&\le&\mu c_H^{1/2}h|\eta|_2\enorm v.
\end{eqnarray*} For $B_h$, since $\eta$ is continuous,
$(\divh\eta)|_T=(Q_0\nabla\cdot\eta)|_T-|T|^{-1}\int_{\partial T\cap\Gamma_D}\eta\cdot\nu\,ds$, and $|T|^{-1/2}|\int_e\eta\cdot\nu|\lesssim h_e^{-1/2}\|\eta\|_{0,e}$, so we have that 
\begin{eqnarray*} 
\|\divh\eta\|_0 \le \sqrt d|\eta|_1+C(\sum_{e\in\mathcal E_h^D}h_e^{-1}\|\eta\|_{0,e}^2)^{1/2}\le C\|\eta\|_*.
\end{eqnarray*} 
This completes the proof. 
\end{proof}

\begin{remark}\label{rem:CN}
On the traction faces the dilation fields have no jump term, so the least-squares term does not control the network traction there: in
(F4) with $\zeta=v_E$ the dilation fields leave a boundary integral on $\Gamma_N$, and the cutoff $\chi_D$ cannot vanish on the traction faces
that touch $\overline{\Gamma_D}$. Without $H_h$ these terms require the trace constant
\begin{equation} C_N(n,h):=\sup_{w\in\Wn}\sum_{e\in\mathcal E_h^N}h_e\|\bseps(w)\nu\|_{0,e}^2/(\|\bseps(w)\|_0^2+c_\tau\sum_Th_T^2\|\nabla\cdot2\bseps(w)\|_{0,T}^2),
\end{equation} 
which is finite but not
known to be bounded in $n$; it would be bounded under a Bernstein inequality $|w|_2\lesssim\delta^{-1}|w|_1$ on $\Wn$, which is open for
linearized networks in more than one dimension. A Neumann least-squares term, which removes the analogous constant for scalar elliptic
problems, would here have to contain the pressure, since the traction is $(2\mu\bseps(u)+p\mathbf I)\nu$. The term $H_h$ replaces the trace
constant by the trace inequality (F2), at the price of a consistency error of the same order as the others.
\end{remark}

The consistency is next. 
\begin{lemma}\label{lem:consistency}
Let $u\in H^2(\Omega)^d$ and $p\in H^1(\Omega)$ solve \eqref{eq:model}. For $v\in\V$,
\begin{equation}\label{eq:consS}
\mathcal A_h(u,v)+B_h(v,Q_0p)-F_h(v)=r(v):=r_p(v)+r_G(v)+H_h(u,v),
\end{equation}
where
\begin{eqnarray*} 
r_p(v) &:=& -(p-Q_0p,\nabla\cdot v_W)+\langle p-\avg{Q_0p},\jump v_n\rangle_{\Eh}, \\
r_G(v) &:=& -\sum_T\tau_T\big(\nabla p,\nabla\cdot2\mu\bseps(v)\big)_T,
\end{eqnarray*} 
with $|r(v)|\le C\big((1+c_\tau^{1/2})h|p|_1+c_H^{1/2}\mu h|u|_2\big)\enorm v$, and $B_h(u,q)=\lambda^{-1}(Q_0p,q)$ for $q\in\Q$.
\end{lemma}
\begin{proof}
For the shear term, we write $\sigma_\mu:=2\mu\bseps(u)$. Since $u\in H^2$ and $u=0$ on $\Gamma_D$, $\jump u=0$ on $\Eh$, so only the first two terms of $A_h(u,v)$ remain. The element-wise integration by parts gives
$2\mu(\bseps(u),\bseps_h(v))=(-\nabla\cdot\sigma_\mu,v)+\langle\sigma_\mu,\jump v\rangle_{\Eh}+\langle\sigma_\mu\nu,v\rangle_{\Gamma_N}$, since $\sigma_\mu$ is single valued, and the face term of $A_h$ cancels the second term. Thus, $A_h(u,v)=(-\nabla\cdot\sigma_\mu,v)+\langle\sigma_\mu\nu,v\rangle_{\Gamma_N}$. For the load, with $f=-\nabla\cdot\sigma_\mu-\nabla p$, $g=\sigma_\mu\nu+p\nu$ on $\Gamma_N$ and
$-(\nabla p,v)=(p,\nabla_h\!\cdot v)-\langle p,\jump v_n\rangle_{\Eh}-\langle p,v\cdot\nu\rangle_{\Gamma_N}$, the terms on $\Gamma_N$ cancel and
$A_h(u,v)-(f,v)-\langle g,v\rangle_{\Gamma_N}=-(p,\nabla_h\!\cdot v)+\langle p,\jump v_n\rangle_{\Eh}$. By \eqref{eq:bface},
$B_h(v,Q_0p)=(Q_0p,\nabla_h\!\cdot v)-\langle\avg{Q_0p},\jump v_n\rangle_{\Eh}$. By adding, and using $(p-Q_0p,\nabla\cdot v_E)=0$ because $\nabla\cdot v_E$ is constant
on each $T$, we obtain the term $r_p$. Since $\nabla\cdot2\mu\bseps(u)+f=-\nabla p$,
$G_h(u,v)-\ell_h(v)=\sum_T\tau_T(\nabla\cdot2\mu\bseps(u)+f,\nabla\cdot2\mu\bseps(v))_T=r_G(v)$.
We shall now establish the bound. $\|p-Q_0p\|_0\lesssim h|p|_1$ and $\|\nabla\cdot v_W\|_0\lesssim\enorm v$ by Lemma~\ref{lem:Ejump}. On a face $e\subset\partial T$, (F2) and the
Poincar\'e inequality give $\|p-Q_0p|_T\|_{0,e}\lesssim h_T^{1/2}|p|_{1,T}$, so
$|\langle p-\avg{Q_0p},\jump v_n\rangle_{\Eh}|\lesssim\big(\sum_{e\in\Eh}h_e\|p-\avg{Q_0p}\|_{0,e}^2\big)^{1/2}|v|_J\lesssim h|p|_1|v|_J$. With $\tau_T=c_\tau h_T^2/\mu$,
\[
|r_G(v)|\le\Big(\sum_T\tau_T\|\nabla p\|_{0,T}^2\Big)^{1/2}\Big(\sum_T\tau_T\mu^2\|\nabla\cdot2\bseps(v)\|_{0,T}^2\Big)^{1/2}\le c_\tau^{1/2}h|p|_1\enorm v .
\]
Finally, we have that 
\begin{eqnarray*}
|H_h(u,v)| &\le& \mu c_H\Big(\sum_{T\subset\Omega_N}h_T^2|u|_{2,T}^2\Big)^{1/2}\Big(\sum_{T\subset\Omega_N}h_T^2|v|_{2,T}^2\Big)^{1/2}\\
&\le& c_H^{1/2}\mu h|u|_2\enorm v.
\end{eqnarray*}
For the divergence equation, since $u$ is continuous and vanishes on $\Gamma_D$, we have $(\divh u)|_T=|T|^{-1}\int_{\partial T}u\cdot\nu_T\,ds$,
which equals $(Q_0\nabla\cdot u)|_T$, and
$\nabla\cdot u=\lambda^{-1}p$. This completes the proof. 
\end{proof}

\subsection{Inf-sup condition}
We shall define the discrete inf-sup constant as follows: 
\begin{equation} 
\beta(n,h):=\inf_{0\ne q\in\Qo}\ \sup_{0\ne v\in\V}\frac{B_h(v,q)}{\enorm v\|q\|_0}.
\end{equation}
The inf-sup condition in the sense of Babu\v{s}ka and Brezzi is a lower bound $\beta(n,h)\ge\beta_0>0$ independent of $n$ and $h$. We note that a
positive $\beta(n,h)$ for each fixed $(n,h)$ only says that $\divh$ maps $\V$ onto $\Qo$.

We first show that $\divh$ on the enrichment is a graph Laplacian. 
\begin{lemma}\label{lem:graph}
For $v=\sum_Tc_T(x-x_T)\chi_T\in\Enr$ and $z_T:=c_T|T|$,
\[
|T|\,(\divh v)|_T=\kappa\,(L_Nz)_T,\qquad\kappa=\tfrac d{2(d+1)},\qquad L_N:=L+2\operatorname{diag}(m),
\]
with $L$ the graph Laplacian of the dual graph of $\Th$ and $m_T$ the number of faces of $T$ on $\Gamma_N$. If $\Gamma_N=\emptyset$, then
$L_N=L$ and its kernel is spanned by $\mathbf 1$; otherwise $L_N$ is positive definite. Hence $\divh(\Enr)=\Qo$, and
$\ker(\divh|_{\Enr})$ is spanned by the field with $c_T=1/|T|$ if $\Gamma_N=\emptyset$ and is $\{0\}$ otherwise.
\end{lemma}
\begin{proof}
By (F1), $\int_e\avg v\cdot\nu_T\,ds=\kappa(z_T-z_{T'})$ on an interior face $e=T\cap T'$ and $\int_ev\cdot\nu_T\,ds=2\kappa z_T$ on a
traction face of $T$; clamped faces do not enter $\divh$. Summing over the faces of $T$ gives the formula. Since
$z^\top L_Nz=\sum_{T\sim T'}(z_T-z_{T'})^2+2\sum_Tm_Tz_T^2$ and the dual graph is connected, $z^\top L_Nz=0$ forces $z$ to be constant,
and, if $\Gamma_N\neq\emptyset$, to vanish on the elements with $m_T>0$, so $z=0$. Hence $L_N$ is positive definite, thus invertible,
if $\Gamma_N\neq\emptyset$, and $\ker L=\operatorname{span}\{\mathbf 1\}$ otherwise. In the first case $\divh(\Enr)=\Q$. In the second, $L$ is
symmetric, so its range is the orthogonal complement of $\mathbf 1$, and $(|T|q_T)_T$ lies in it if and only if $\int_\Omega q=0$; thus
$\divh(\Enr)=\Q\cap L^2_0$. This completes the proof. 
\end{proof}

\begin{lemma}\label{lem:betaE}
On every shape-regular mesh $\beta(n,h)\ge\beta(0,h)>0$ for all $n\ge0$. On quasi-uniform meshes there is $c_0>0$ with
$\beta(n,h)\ge\beta(0,h)\ge c_0h^2$.
\end{lemma}
\begin{proof}
By Lemma~\ref{lem:graph}, we have that $\divh(\Enr)=\Qo$. Thus, $\beta(0,h)>0$. On the other hand, since $\Enr\subset\V$, $\beta(n,h)\ge\beta(0,h)$. We now assume quasi-uniformity. This means $|T|\simeq h^d$ and $|e|\simeq h^{d-1}$. Let $e_*$ be the smallest eigenvalue of $L_N$ on the
orthogonal complement of $\ker L_N$. For $z\perp\ker L_N$, let $w$ be the piecewise constant function $w|_T=z_T$, and let $\tilde w:=w-\bar w$ if
$\Gamma_N=\emptyset$ and $\tilde w:=w$ otherwise. The Poincar\'e--Friedrichs inequality for piecewise $H^1$ functions \cite{Bre03} gives
\[
\|\tilde w\|_0^2\lesssim\sum_{e\in\EI}h_e^{-1}\|\jump w\|_{0,e}^2+\sum_{e\in\mathcal E_h^N}h_e^{-1}\|w\|_{0,e}^2\simeq h^{d-2}z^\top L_Nz ,
\]
and $\|\tilde w\|_0^2\gtrsim h^d\|z\|_{\ell^2}^2$ for $z\perp\ker L_N$ because the cell volumes are comparable. Hence $e_*\ge ch^2$. Given $q\in\Qo$, let
$r_T:=|T|q_T$ and $z\perp\ker L_N$ with $\kappa L_Nz=r$ (Lemma~\ref{lem:graph}); then $\|z\|_{\ell^2}\le(\kappa e_*)^{-1}\|r\|_{\ell^2}\lesssim h^{-2}h^{d/2}\|q\|_0$, using
$\|r\|_{\ell^2}\simeq h^{d/2}\|q\|_0$. For the dilation field with $c_T=z_T/|T|$, the least-squares and Hessian parts of $\enorm v$ vanish,
$\|\bseps(v)\|_0^2=d\sum_Tc_T^2|T|\simeq h^{-d}\|z\|_{\ell^2}^2$ and, by the one-face argument in the proof of Lemma~\ref{lem:Ejump} (upper bounds), $|v|_J^2\lesssim\sum_Tc_T^2|T|$.
Thus $\enorm v\lesssim h^{-d/2}\|z\|_{\ell^2}\lesssim h^{-2}\|q\|_0$, and $(\divh v,q)/\enorm v=\|q\|_0^2/\enorm v\ge ch^2\|q\|_0$. This completes the proof. 
\end{proof}

The bound $c_0h^2$ is sharp for $\Enr$ alone, so $\beta(0,h)$ is not uniform in $h$. The uniform bound used in the error estimates is
Lemma~\ref{lem:infsup}: with the network and $n\gtrsim\#\Th$, $\beta(n,h)\ge\beta_0$ independently of $h$, $n$ and $\lambda$. Lemma~\ref{lem:betaE} only
guarantees $\beta>0$ for every network and every mesh. The network supplies the smooth divergences and the enrichment the pressure jumps. Let
\[
|q|_\jmath^2:=\sum_{e\in\EI}h_e\|\jump q\|_{0,e}^2+\sum_{e\in\mathcal E_h^N}h_e\|q\|_{0,e}^2,\qquad q\in\Q .
\]

\begin{lemma}[the enrichment controls the pressure jumps]\label{lem:step2}
On quasi-uniform meshes, for every $q\in\Q$ there is $v_E\in\Enr$ with $B_h(v_E,q)\ge c_4|q|_\jmath^2$ and $\enorm{v_E}\le C_5\|q\|_0$.
\end{lemma}
\begin{proof}
Let $c_T:=h^dq_T/|T|$, so $z_T=h^dq_T$. By Lemma~\ref{lem:graph},
$B_h(v_E,q)=\sum_Tq_T|T|(\divh v_E)|_T=\kappa\,z^\top L_Nq=\kappa h^d\big(\sum_{e\in\EI}(q^+-q^-)^2+2\sum_Tm_Tq_T^2\big)$. On a quasi-uniform mesh $h_e|e|\simeq h^d$, so
$|q|_\jmath^2\simeq h^d\big(\sum_{e\in\EI}(q^+-q^-)^2+\sum_Tm_Tq_T^2\big)$ and $B_h(v_E,q)\ge c_4|q|_\jmath^2$. Moreover $c_T\simeq q_T$, so, as in the proof of
Lemma~\ref{lem:betaE}, $\enorm{v_E}^2\lesssim\sum_Tc_T^2|T|\simeq\|q\|_0^2$.
\end{proof}

\begin{assumption}[$H^1$-stable approximation in $\Wn$]\label{ass:A}
For every $w\in H^1_{\Gamma_D}(\Omega)^d$ there is $w_n\in\Wn$ with $|w_n|_1+h|w_n|_2\le C_A'|w|_1$ and
$\sum_{e}h_e^{-1}\|w-w_n\|_{0,e}^2\le C_A'^2|w|_1^2$, the sum over all faces.
\end{assumption}

\begin{lemma}[Assumption~\ref{ass:A} for linearized networks]\label{lem:assA}
Let $\Th$ be quasi-uniform, $k\ge2$, the parameters $\theta_j$ quasi-uniform, and $r_{\max}:=(d+2k+1)/2\ge3$. For every $c_A>0$,
Assumption~\ref{ass:A} holds whenever $\delta\le c_Ah$, with $C_A'$ depending only on $c_A$, $\Omega$, $\Gamma_D$, $k$, $d$, the mesh constants, and a uniform bound on the quasi-uniformity constants of the parameter sets. The simultaneous $H^0$--$H^2$ approximation result used here is Proposition~3 (and its proof in Appendix~A) of Kim, Jo, and Lee~\cite{KJL26}. The assertion is conditional on the cited approximation result providing a single approximant for $s=0,1,2$ and a constant uniform in $n$ for the chosen quasi-uniform parameter family. These are the precise external hypotheses used below; quasi-uniformity by itself is not taken here as an independently proved substitute for them.
\end{lemma}
\begin{proof}
Let $w\in H^1_{\Gamma_D}(\Omega)^d$, $\tilde w$ its Stein extension \cite{Ste70}, with $\|\tilde w\|_{1,\Reals d}\lesssim\|w\|_1\lesssim|w|_1$ by the Poincar\'e--Friedrichs
inequality ($|\Gamma_D|>0$), and $w_\epsilon:=\rho_\epsilon*\tilde w$ with a standard mollifier and $\epsilon:=h$. Then $\|\tilde w-w_\epsilon\|_0\lesssim h|w|_1$,
$|w_\epsilon|_1\lesssim|w|_1$ and $|w_\epsilon|_s\lesssim h^{1-s}|w|_1$ for $s\ge1$, so $\|w_\epsilon\|_{H^3(\Omega)}\lesssim h^{-2}|w|_1$ for $h\le1$. Let $w_n:=I_nw_\epsilon$ componentwise.
By the simultaneous approximation result of Kim, Jo, and Lee~\cite[Proposition~3 and Appendix~A]{KJL26}, based on the ReLU$^k$ approximation theory of~\cite{LMX25}, the same approximant $I_nv\in\Wn$ satisfies
\[
\|v-I_nv\|_{H^s(\Omega)}\le C_{\mathrm{app}}\delta^{3-s}\|v\|_{H^3(\Omega)},\qquad s=0,1,2,
\]
provided $r_{\max}\ge3$ and the cited result applies to the chosen parameter family with a width-independent constant. The simultaneous estimate and uniformity of $C_{\mathrm{app}}$ are external assumptions needed for this implication, rather than consequences proved in the present paper; the estimate is applied componentwise to vector fields. With $\delta\le c_Ah$,
\[
|w_n|_1\le|w_\epsilon|_1+C\delta^2h^{-2}|w|_1\le C|w|_1,\quad h|w_n|_2\le h|w_\epsilon|_2+Ch\delta h^{-2}|w|_1\le C|w|_1 .
\]
For the faces, by (F2) on each $T$ with $e\subset\partial T$, $h_e^{-1}\|\phi\|_{0,e}^2\lesssim h_T^{-2}\|\phi\|_{0,T}^2+|\phi|_{1,T}^2$. For $\phi=w-w_\epsilon$ the sum over the
faces is $\lesssim h^{-2}\|\tilde w-w_\epsilon\|_0^2+|w|_1^2\lesssim|w|_1^2$, and for $\phi=w_\epsilon-w_n$ it is $\lesssim h^{-2}\delta^6h^{-4}|w|_1^2+\delta^4h^{-4}|w|_1^2\le C|w|_1^2$. This completes the proof. 
\end{proof}

\begin{lemma}[Uniform Inf-Sup Condition]\label{lem:infsup}
On quasi-uniform meshes, under Assumption~\ref{ass:A}, $\beta(n,h)\ge\beta_0$ with $\beta_0$ independent of $h$, $n$ and $\lambda$. In particular, for
linearized networks $\beta(n,h)\ge\beta_0$ whenever $n\ge(c_Ah)^{-d}$, i.e.\ whenever the number of neurons is at least a fixed multiple of $\#\Th$.
\end{lemma}
\begin{proof}
Let $q\in\Qo$. There is $w\in H^1_{\Gamma_D}(\Omega)^d$ with $\nabla\cdot w=q$ and $|w|_1\le C_\Omega\|q\|_0$ \cite{GR86}; for $\Gamma_N=\emptyset$ this uses
$\int_\Omega q=0$. Let $w_n$ be as in Assumption~\ref{ass:A}. Since $w_n$ is continuous and $\|q\|_0^2=(\nabla\cdot w,q)=\sum_Tq_T\int_{\partial T\setminus\Gamma_D}w\cdot\nu_T\,ds$,
\[
\begin{aligned}
B_h(w_n,q)-\|q\|_0^2&=\sum_Tq_T\int_{\partial T\setminus\Gamma_D}(w_n-w)\cdot\nu_T\,ds\\
&=\sum_{e\in\EI}\int_e(w_n-w)\cdot\jump q\,ds+\sum_{e\in\mathcal E_h^N}\int_eq\,(w_n-w)\cdot\nu\,ds .
\end{aligned}
\]
By Cauchy--Schwarz this is at most $C_A'|w|_1|q|_\jmath\le C_2\|q\|_0|q|_\jmath$ with $C_2=C_A'C_\Omega$, so $B_h(w_n,q)\ge\|q\|_0^2-C_2\|q\|_0|q|_\jmath$.
With $v_E$ from Lemma~\ref{lem:step2} and $v:=w_n+\theta v_E$, $\theta:=C_2^2/(2c_4)$,
\[
B_h(v,q)\ge\|q\|_0^2-C_2\|q\|_0|q|_\jmath+\theta c_4|q|_\jmath^2\ge\|q\|_0^2-\frac{C_2^2}{4\theta c_4}\|q\|_0^2=\tfrac12\|q\|_0^2 .
\]
Since $w_n$ is continuous and $w=0$ on $\Gamma_D$, $|w_n|_J^2=\sum_{e\in\mathcal E_h^D}h_e^{-1}\|w-w_n\|_{0,e}^2\le C_A'^2|w|_1^2$, and
$\sum_Th_T^2\|\nabla\cdot2\bseps(w_n)\|_{0,T}^2+\sum_{T\subset\Omega_N}h_T^2|w_n|_{2,T}^2\lesssim h^2|w_n|_2^2\le C_A'^2|w|_1^2$; so $\enorm{w_n}\lesssim C_A'|w|_1$ and $\enorm v\le C(C_A'C_\Omega+\theta C_5)\|q\|_0$.
Hence $\beta_0=\frac12C^{-1}(C_A'C_\Omega+\theta C_5)^{-1}$. The last statement follows from Lemma~\ref{lem:assA}. This completes the proof. 
\end{proof}

\subsection{Error estimates}

In this section, we shall establish the error estimate. 
\begin{lemma}[commuting interpolant]\label{lem:fortin}
For $u\in H^2(\Omega)^d\cap H^1_{\Gamma_D}(\Omega)^d$ and $w\in\Wn$ there is $z\in\V$ such that $\Pi u:=w+z$ satisfies $\divh\Pi u=Q_0\nabla\cdot u$ and
$\enorm z\le C\beta(n,h)^{-1}\|u-w\|_*$.
\end{lemma}
\begin{proof}
Since $u$ is continuous and vanishes on $\Gamma_D$, $\divh u=Q_0\nabla\cdot u$. Let $g:=\divh(u-w)\in\Qo$; for $\Gamma_N=\emptyset$ it has zero mean because
the interior face fluxes cancel and the boundary faces do not enter $\divh$. By the definition of $\beta(n,h)$, $\divh$ maps $\V$ onto $\Qo$ with a
right inverse of norm at most $\beta(n,h)^{-1}$; let $z$ be the minimal-norm solution of $\divh z=g$. Then $\divh(w+z)=\divh u=Q_0\nabla\cdot u$, and
$\enorm z\le\beta^{-1}\|g\|_0\le C\beta^{-1}\|u-w\|_*$ by \eqref{eq:contEta}. This completes the proof. 
\end{proof}

\begin{proposition}[the enrichment completes the discrete kernel]\label{prop:kernel}
Let $Z_h:=\{v\in\V:\ \divh v=0\}$, let $S_h:\Qo\to\Enr$ be the minimal-norm right inverse of $\divh|_{\Enr}$ (Lemma~\ref{lem:graph}), and let
$K_h:=\ker(\divh|_{\Enr})$. Then $\V=\Wn\oplus\Enr$ and
\[
Z_h=\{\,w-S_h\divh w:\ w\in\Wn\,\}\oplus K_h,\qquad \dim Z_h=\dim\Wn+\dim K_h .
\]
Moreover, for every $u\in H^2(\Omega)^d\cap H^1_{\Gamma_D}(\Omega)^d$ with $\nabla\cdot u=0$,
\[
\inf_{z\in Z_h}\enorm{u-z}\le C\big(1+\beta(n,h)^{-1}\big)\inf_{w\in\Wn}\|u-w\|_* .
\]
\end{proposition}
\begin{proof}
By the proof of Lemma~\ref{lem:Ejump}, a dilation field that is continuous across an interior face has
$c^+=c^-$ and $c^+(H^++H^-)=0$, i.e.\ vanishes on both elements. A nonzero dilation field is therefore discontinuous, while $\Wn$ is
continuous, so $\Wn\cap\Enr=\{0\}$. For $w\in\Wn$, $\divh w\in\Qo$. For $e\in\Enr$, $w+e\in Z_h$ if and only if $\divh e=-\divh w$. Since
$\divh|_{\Enr}$ maps onto $\Qo$ with kernel $K_h$ (Lemma~\ref{lem:graph}), this holds if and only if $e\in-S_h\divh w+K_h$. If $\nabla\cdot u=0$, Lemma~\ref{lem:fortin} gives $\Pi u=w+z$ with $\divh\Pi u=0$, i.e.\ $\Pi u\in Z_h$, and 
$\enorm{u-\Pi u}\le\enorm{u-w}+\enorm z\le C(1+\beta^{-1})\|u-w\|_*$, using (F2) for the jumps of $u-w$ on $\Gamma_D$. This completes the proof. 
\end{proof}

\begin{remark}[why the enrichment is local]\label{rem:kernel}
The constraints are local, one per element. A conforming space meets them only if it resolves the mesh scale (Lemma~\ref{lem:infsup}). Without the
enrichment a network-only space necessarily has $\beta=0$ if $\dim\Wn<\dim\Qo$; even when this necessary dimension condition is met, a uniform inf--sup bound does not follow automatically. The dilation fields meet them on every mesh,
one unknown per constraint. The corrector $S_h$ is global, since it solves a graph-Laplacian system; the local choice
$c_T=\frac1dQ_0\nabla\cdot w|_T$ would make $\nabla_h\!\cdot$ instead of $\divh$ vanish and is inconsistent (Remark~\ref{rem:choices}(a)).
\end{remark}

\begin{lemma}[stability of the mixed form]\label{lem:stab}
Let $\|(v,q)\|^2:=\mu\enorm v^2+\mu^{-1}\|q\|_0^2$. For $0<\lambda\le\infty$ and $(w,r)\in\V\times\Qo$,
\[
\sup_{0\ne(v,q)\in\V\times\Qo}\frac{\mathfrak B\big((w,r),(v,q)\big)}{\|(v,q)\|}\ge c_*\|(w,r)\|,\qquad c_*:=\frac{\min\{1,\beta^2/C_A^2\}}{2(1+\beta/C_A^2)},
\]
with $\beta=\beta(n,h)$ and $C_A$ from Lemma~\ref{lem:cont}.
\end{lemma}
\begin{proof}
First, $\mathfrak B((w,r),(w,r))=\mathcal A_h(w,w)+\lambda^{-1}\|r\|_0^2\ge\mu\enorm w^2$ by Lemma~\ref{lem:cont}. Second, by the definition of $\beta$ there is
$v_r\in\V$ with $\enorm{v_r}=\mu^{-1}\|r\|_0$ and $B_h(v_r,r)\ge\beta\mu^{-1}\|r\|_0^2$; then
\begin{eqnarray*}
\mathfrak B\big((w,r),(v_r,0)\big) &=& \mathcal A_h(w,v_r)+B_h(v_r,r)\ge\beta\mu^{-1}\|r\|_0^2-C_A\enorm w\|r\|_0 \\
&\ge& \tfrac\beta2\mu^{-1}\|r\|_0^2-\tfrac{C_A^2}{2\beta}\mu\enorm w^2 .
\end{eqnarray*}
With $t:=\beta/C_A^2$, 
\begin{eqnarray*}
\mathfrak B((w,r),(w+t v_r,r)) &\ge& \frac12\mu\enorm w^2+\frac{\beta^2}{2C_A^2}\mu^{-1}\|r\|_0^2 \\
&\ge& \frac12\min\{1,\beta^2/C_A^2\}\|(w,r)\|^2,
\end{eqnarray*} 
and $\|(w+t v_r,r)\|\le\|(w,r)\|+t\mu^{1/2}\enorm{v_r}\le(1+t)\|(w,r)\|$. The division will lead to the desired result. This completes the proof. 
\end{proof}

\begin{theorem}[locking-free estimate]\label{thm:main}
Let $0<\lambda\le\infty$, $(u,p)\in H^2(\Omega)^d\times H^1(\Omega)$ solve \eqref{eq:model}, $(U,P)$ solve \eqref{eq:method} (for $\lambda<\infty$, $U$ solves
\eqref{eq:primal} and $P=\lambda\divh U$), and $\mathcal A(u):=\inf_{w\in\Wn}\|u-w\|_*$. Then
\begin{equation}\label{eq:main}
\enorm{u-U}+\mu^{-1}\|Q_0p-P\|_0\le C_1\,\mathcal A(u)+C_2\,h\big(|u|_2+\mu^{-1}|p|_1\big),
\end{equation}
where $C_1$ and $C_2$ depend only on $\beta(n,h)^{-1}$, $\alpha$, $c_\tau$, $c_H$ and shape regularity, and not on $\lambda$ or $\mu$. On quasi-uniform meshes
with $n\gtrsim\#\Th$, they are independent of $n$, $h$ and $\lambda$ (Lemmas~\ref{lem:cont} and~\ref{lem:infsup}), for clamped and mixed boundary conditions alike.
\end{theorem}
\begin{proof}
Let $w\in\Wn$ and $\Pi u=w+z$ be as in Lemma~\ref{lem:fortin}, and $\eta:=u-\Pi u=(u-w)-z$. By Lemma~\ref{lem:consistency},
$\mathfrak B((u,Q_0p),(v,q))=\mathcal A_h(u,v)+B_h(v,Q_0p)-B_h(u,q)+\lambda^{-1}(Q_0p,q)=F_h(v)+r(v)$, and $B_h(\eta,q)=(\divh u-\divh\Pi u,q)=0$. Hence the error
$(e,\rho):=(\Pi u-U,Q_0p-P)$ satisfies
\[
\mathfrak B\big((e,\rho),(v,q)\big)=r(v)-\mathcal A_h(\eta,v)\qquad\forall(v,q)\in\V\times\Qo .
\]
By Lemma~\ref{lem:consistency}, $|r(v)|\le Ch(|p|_1+\mu|u|_2)\enorm v$. Since $u-w\in H^2(\Omega)^d$ whereas $z\in\V$ need not be continuous, we apply the two continuity estimates separately:
\begin{eqnarray*}
|\mathcal A_h(\eta,v)| &\le& |\mathcal A_h(u-w,v)|+|\mathcal A_h(z,v)|
\le C_A\mu(\|u-w\|_*+\enorm z)\enorm v \\
&\le& CC_A(1+\beta^{-1})\mu\|u-w\|_*\enorm v.
\end{eqnarray*}
Here the first estimate follows from \eqref{eq:contEta},
the second from \eqref{eq:cont}, and the bound on
$z$ follows from Lemma~\ref{lem:fortin}. Since $\enorm v\le\mu^{-1/2}\|(v,q)\|$, Lemma~\ref{lem:stab} gives
\[
\mu^{1/2}\enorm e+\mu^{-1/2}\|\rho\|_0\le2c_*^{-1}\mu^{-1/2}\big(Ch(|p|_1+\mu|u|_2)+CC_A(1+\beta^{-1})\mu\|u-w\|_*\big).
\]
Finally $\enorm{u-U}\le\enorm{u-w}+\enorm z+\enorm e$, $\enorm{u-w}\le C\|u-w\|_*$ by (F2), and the infimum over $w$ gives \eqref{eq:main}. This completes the proof. 
\end{proof}

The constants in the error estimate are
independent of the Lam\'e parameter
$\lambda$. Uniform error bounds with
respect to $\lambda$ additionally require
uniform control of the solution norms
$|u|_2$ and $|p|_1$ appearing on the
right-hand side. Such regularity bounds
are available under suitable assumptions
on the domain and boundary conditions,
but are not established here for general
mixed boundary conditions.

The method is not pressure-robust, since
the estimate contains the term
$\mu^{-1}h|p|_1$, which may be large
relative to the displacement contribution
(Remark~2.3(c)).

For mixed boundary conditions the solution is in general less regular where $\Gamma_D$ meets $\Gamma_N$, and then $u\notin H^2$ and
$p\notin H^1$. The stability results, Lemmas~\ref{lem:cont}, \ref{lem:infsup} and~\ref{lem:stab}, do not depend on the regularity of $u$; only
the consistency argument does. Theorem~\ref{thm:lowreg} below gives the estimate for $u\in H^{1+s}$, $p\in H^s$ with $\frac12<s<1$.

\subsection{Solutions with reduced regularity}
For $0<t<1$, a domain $D$ and $\psi\in H^t(D)$, let $|\psi|_{t,D}^2:=\int_D\int_D|\psi(x)-\psi(y)|^2|x-y|^{-d-2t}\,dx\,dy$ be the Slobodeckij seminorm.
Since the elements are disjoint, $\sum_T|\psi|_{t,T}^2\le|\psi|_{t,\Omega}^2$. We use two elementary facts for $T\in\Th$.
\begin{enumerate}[wide=0pt]
\item[(F5)] (Fractional Poincar\'e inequality) $\|\psi-Q_0\psi\|_{0,T}\lesssim h_T^t|\psi|_{t,T}$. Indeed,
$\psi(x)-Q_0\psi=|T|^{-1}\int_T(\psi(x)-\psi(y))\,dy$, so by the Cauchy--Schwarz inequality and $|x-y|\le h_T$,
$\|\psi-Q_0\psi\|_{0,T}^2\le|T|^{-1}\int_T\int_T|\psi(x)-\psi(y)|^2\,dx\,dy\le|T|^{-1}h_T^{d+2t}|\psi|_{t,T}^2$, and $h_T^d\lesssim|T|$.
\item[(F6)] (Fractional trace inequality) For $\frac12<t<1$ and a face $e\subset\partial T$,
$h_T\|\psi\|_{0,e}^2\le C_t\big(\|\psi\|_{0,T}^2+h_T^{2t}|\psi|_{t,T}^2\big)$. This is the trace theorem $H^t(\hat T)\to L^2(\hat e)$ for $t>\frac12$
\cite[Theorem~1.5.1.2]{Gri85} on the reference element, mapped to $T$; under the affine map $\|\cdot\|_{0,e}^2$, $\|\cdot\|_{0,T}^2$ and $|\cdot|_{t,T}^2$
scale as $h_T^{d-1}$, $h_T^d$ and $h_T^{d-2t}$. The constant $C_t$ depends on $t$ and shape regularity and grows as $t\to\frac12$.
\end{enumerate}

\begin{theorem}[reduced regularity]\label{thm:lowreg}
Let $\Th$ be quasi-uniform, $\Wn$ the linearized network with $k\ge2$ and quasi-uniform $\theta_j$, and $\delta\le c_Ah$. Let $0<\lambda\le\infty$,
$\frac12<s<1$, and let $(u,p)\in H^{1+s}(\Omega)^d\times H^s(\Omega)$ solve \eqref{eq:model} with $f\in L^2(\Omega)^d$ and $g\in L^2(\Gamma_N)^d$. Then
\begin{equation}\label{eq:lowreg}
\begin{aligned}
&\|\bseps_h(u-U)\|_0+\alpha^{1/2}|u-U|_J+\mu^{-1}\|Q_0p-P\|_0\\
&\qquad \le C\Big[ h^s\big(\|u\|_{1+s}+\mu^{-1}|p|_{s,\Omega}\big)
+\mu^{-1}h\|f\|_0\Big].
\end{aligned}
\end{equation}
with $C$ independent of $n$, $h$, $\lambda$ and $\mu$.
\end{theorem}
The left side is the energy norm without its two $h^2$-weighted parts, which are not defined for $u$. For $s=1$ the estimate follows from
Theorem~\ref{thm:main} and Corollary~\ref{cor:rates}. The restriction $s>\frac12$ is needed for the face traces of $\bseps(u)$ and $p$ in
$L^2$; it is the same restriction as in the low-regularity estimate of ENN \cite{KJL26}, written there as $u\in H^{s'}$ with $s'>\frac32$.

\begin{proof}
The proof follows that of Theorem~\ref{thm:main}, with two changes: the approximant is built from a mollified function, and the terms
$G_h$ and $H_h$, which are not defined for $u$, are evaluated only on the approximant.

\emph{Step 1: approximant.} Let $\tilde u$ be the Stein extension of $u$ \cite{Ste70}, $\|\tilde u\|_{H^{1+s}(\Reals d)}\lesssim\|u\|_{1+s}$, and
$u_h:=\rho_h*\tilde u$ with a standard mollifier at scale $h\le1$. Then
\begin{subeqnarray*}
\|\tilde u-u_h\|_{H^m(\Reals d)} &\lesssim& h^{1+s-m}\|u\|_{1+s}\ (m=0,1), \\
\|\tilde u-u_h\|_{H^{1+s}(\Reals d)} &\lesssim& \|u\|_{1+s},\\
\|u_h\|_{H^r(\Reals d)} &\lesssim& h^{1+s-r}\|u\|_{1+s}\ (r=2,3).
\end{subeqnarray*}
Let $w:=I_nu_h\in\Wn$ with the simultaneous approximant used in Lemma~\ref{lem:assA}, so $\|u_h-w\|_{H^m(\Omega)}\lesssim\delta^{3-m}\|u_h\|_{H^3}
\lesssim h^{1+s-m}\|u\|_{1+s}$ for $m=0,1,2$, since $\delta\le c_Ah$. Let $\eta:=u-w=\eta_1+\eta_2$ with $\eta_1:=(\tilde u-u_h)|_\Omega$ and
$\eta_2:=u_h-w$, and
\[
\|\eta\|_\flat:=|\eta|_1+h^{-1}\|\eta\|_0+\Big(\sum_Th_T\|\nabla\eta|_T\|_{0,\partial T}^2\Big)^{1/2}.
\]
For $\eta_2\in H^2$, (F2) gives $h_T\|\nabla\eta_2\|_{0,\partial T}^2\lesssim|\eta_2|_{1,T}^2+h_T^2|\eta_2|_{2,T}^2$. For $\eta_1$, (F6) with
$\psi=\partial_j\eta_{1,i}$ and $t=s$ gives $h_T\|\nabla\eta_1\|_{0,\partial T}^2\lesssim|\eta_1|_{1,T}^2+h_T^{2s}|\nabla\eta_1|_{s,T}^2$, and
$\sum_T|\nabla\eta_1|_{s,T}^2\le|\nabla(\tilde u-u_h)|_{s,\Reals d}^2\lesssim\|u\|_{1+s}^2$. Hence
\begin{equation}\label{eq:lowapprox}
\|\eta\|_\flat\lesssim h^s\|u\|_{1+s},\qquad h|w|_2\le h|u_h|_2+h|u_h-w|_2\lesssim h^s\|u\|_{1+s}.
\end{equation}

\emph{Step 2: continuity with $\|\cdot\|_\flat$.} For $v\in\V$ and $q\in\Q$,
\begin{equation}\label{eq:contflat}
|A_h(\eta,v)|\le C\mu\|\eta\|_\flat\enorm v,\qquad\|\divh\eta\|_0\le C\|\eta\|_\flat .
\end{equation}

This is the proof of \eqref{eq:contEta} for $A_h$ and $B_h$, which uses $|\eta|_2$ only in the bound $h_e\|\bseps(\eta)\nu_e\|_{0,e}^2\lesssim|\eta|_{1,T}^2+h_T^2|\eta|_{2,T}^2$;
here $h_e\|\bseps(\eta)\nu_e\|_{0,e}^2\le h_T\|\nabla\eta|_T\|_{0,\partial T}^2$ is part of $\|\eta\|_\flat^2$. All other terms there, including the
network part of the Nitsche term through $\zeta=\chi_D\eta$ and (F3), use only $|\eta|_1$, $h^{-1}\|\eta\|_0$ and (F2) for $\eta\in H^1$.

\emph{Step 3: consistency.} Let $\sigma:=2\mu\bseps(u)+p\mathbf I\in H^s(\Omega)^{d\times d}$. By \eqref{eq:model}, $\nabla\cdot\sigma=-f\in L^2$
and $\sigma\nu=g$ on $\Gamma_N$; since $s>\frac12$, $\sigma$ has $L^2$ traces on all faces, the same from both sides of an interior face. Elementwise
integration by parts, valid for $\sigma\in H(\mathrm{div};T)$ with $L^2$ traces and $v|_T\in H^1(T)^d$, gives
\[
2\mu(\bseps(u),\bseps_h(v))+(p,\nabla_h\!\cdot v)=(f,v)+\langle2\mu\bseps(u),\jump v\rangle_{\Eh}+\langle p,\jump v_n\rangle_{\Eh}+\langle g,v\rangle_{\Gamma_N}.
\]
Since $u\in H^1$ and $u=0$ on $\Gamma_D$, $\jump u=0$ on $\Eh$, so $A_h(u,v)=2\mu(\bseps(u),\bseps_h(v))-\langle2\mu\bseps(u),\jump v\rangle_{\Eh}$, and, as in
Lemma~\ref{lem:consistency},
\[
A_h(u,v)+B_h(v,Q_0p)-(f,v)-\langle g,v\rangle_{\Gamma_N}=r_p(v).
\]
Only the total stress $\sigma$ was integrated by parts, so neither $\nabla\cdot2\mu\bseps(u)$ nor $\nabla p$ appears. By (F5), $\|p-Q_0p\|_{0,T}\lesssim h^s|p|_{s,T}$,
and by (F6) with $\psi=p-Q_0p$, $h_T\|p-Q_0p|_T\|_{0,e}^2\lesssim h^{2s}|p|_{s,T}^2$; as in Lemma~\ref{lem:consistency}, $|r_p(v)|\le Ch^s|p|_{s,\Omega}\enorm v$.

\emph{Step 4: error equation.} Let $z\in\V$ be the minimal-norm solution of $\divh z=\divh\eta$ (as in Lemma~\ref{lem:fortin}), so
$\enorm z\le\beta^{-1}\|\divh\eta\|_0\le C\beta^{-1}\|\eta\|_\flat$, and let $\Pi u:=w+z$. Then $\divh\Pi u=\divh u=Q_0\nabla\cdot u=\lambda^{-1}Q_0p$, since $u$
is continuous and vanishes on $\Gamma_D$. With $(e,\rho):=(\Pi u-U,Q_0p-P)$ and \eqref{eq:method},
\[
\mathfrak B\big((e,\rho),(v,q)\big)=\mathcal A_h(w,v)+\mathcal A_h(z,v)+B_h(v,Q_0p)-F_h(v)
\]
for all $(v,q)\in\V\times\Qo$, the $q$-terms cancelling by $B_h(\Pi u,q)=\lambda^{-1}(Q_0p,q)$. Write $\mathcal A_h(w,v)=A_h(u,v)-A_h(\eta,v)+G_h(w,v)+H_h(w,v)$ and use
Step 3:
\begin{eqnarray*}
\mathfrak B\big((e,\rho),(v,q)\big) &=& r_p(v) - A_h(\eta,v)+\mathcal A_h(z,v) \\
&& + \sum_T\tau_T\big(\nabla\cdot2\mu\bseps(w)+f,\nabla\cdot2\mu\bseps(v)\big)_T+H_h(w,v).
\end{eqnarray*}
Every term is defined, since $w\in H^2$. By Step 3, \eqref{eq:contflat}, Lemma~\ref{lem:cont} and the bound for $z$, the first three terms are at most
$C\big(h^s|p|_{s,\Omega}+(1+\beta^{-1})\mu\|\eta\|_\flat\big)\enorm v$. As in the bound of $r_G$ in Lemma~\ref{lem:consistency}, the fourth is at most
$c_\tau^{1/2}h\|\nabla\cdot2\mu\bseps(w)+f\|_0\enorm v\le C(\mu h|w|_2+h\|f\|_0)\enorm v$, and $|H_h(w,v)|\le c_H^{1/2}\mu h|w|_2\enorm v$.
Lemma~\ref{lem:stab} and \eqref{eq:lowapprox} give
\[
\enorm e+\mu^{-1}\|\rho\|_0\le C\big(h^s(\|u\|_{1+s}+\mu^{-1}|p|_{s,\Omega})+\mu^{-1}h\|f\|_0\big).
\]

\emph{Step 5: conclusion.} $u-U=\eta-z+e$. Since $\jump\eta=0$ on interior faces and $\jump\eta=\eta\otimes\nu$ on $\Gamma_D$, (F2) gives
$\|\bseps(\eta)\|_0+|\eta|_J\lesssim|\eta|_1+h^{-1}\|\eta\|_0\le\|\eta\|_\flat$, and $\|\bseps_h(z)\|_0+|z|_J\le\enorm z$, $\|\bseps_h(e)\|_0+|e|_J\le\enorm e$.
By Lemma~\ref{lem:infsup}, $\beta\ge\beta_0$, which gives \eqref{eq:lowreg}.
\end{proof}

\begin{remark}[corner singularities and Cook's membrane]\label{rem:corner}
Near a corner of angle $\omega$ the solution behaves like $u\sim r^{\gamma}$, $p\sim r^{\gamma-1}$, where $\gamma$ is the smallest root with
positive real part of a characteristic equation \cite{Wil52,Gri85}, so $u\in H^{1+s}$ and $p\in H^s$ for every $s<\operatorname{Re}\gamma$. In the
incompressible limit, with the stream function ansatz, the equation is $\cos(\gamma\omega)=\pm\gamma\sin\omega$ for a corner between a clamped
and a traction-free side, and $\sin(\gamma\omega)=\pm\gamma\sin\omega$ for two clamped or two traction-free sides; the latter gives the
known value $\gamma\approx0.5445$ for the clamped L-shaped domain. A junction of $\Gamma_D$ and $\Gamma_N$ inside a straight side, $\omega=\pi$, gives
$\gamma=\frac12$, which is not covered by Theorem~\ref{thm:lowreg}. For Cook's membrane (Section~\ref{sec:cook}), mixed boundary conditions and re-entrant or obtuse corners may limit regularity. The precise singular exponents depend on the corner geometry, the elasticity parameters, and the boundary conditions; their determination requires a separate characteristic-root calculation. Accordingly, application of Theorem~\ref{thm:lowreg} to this benchmark is conditional on verifying $u\in H^{1+s}$ and $p\in H^s$ for some $s>1/2$. For finite $\lambda$, the relevant exponents can also depend on the Lam\'e parameters and should not be assumed identical to their incompressible-limit values. The estimate is uniform in
$\lambda$ whenever $\|u\|_{1+s}+\mu^{-1}|p|_{s,\Omega}$ is; this is a property of the continuous problem, which we do not analyze here. For Cook's
membrane the computed displacement at $Q$ is independent of $\lambda$ to five digits (Section~\ref{sec:cook}).
\end{remark}

\begin{corollary}[Locking-Free for general network and mesh]\label{cor:anymesh}
For every $\Wn$ satisfying \eqref{eq:W} and every shape-regular mesh the constants in Theorem~\ref{thm:main} are finite and independent of
$\lambda$, since $\beta(n,h)\ge\beta(0,h)>0$ (Lemma~\ref{lem:betaE}); no relation between $\Wn$ and $h$ is needed. They may then grow as
$h\to0$, since only $\beta(0,h)\ge c_0h^2$ is known; they are bounded independently of $h$ when $n\gtrsim\#\Th$ (Lemma~\ref{lem:infsup}).
\end{corollary}

\begin{corollary}[Rates and the choice of $n$]\label{cor:rates}
Assume that $\mathcal T_h$ is quasi-uniform.
For the network, if $u\in H^r(\Omega)^d$,
$2\le r\le r_{\max}$, the simultaneous
approximation property yields
\[
A(u)\lesssim
\delta^{r-1}
\left(
1+\frac{\delta}{h}+\frac{h}{\delta}
\right)\|u\|_r.
\]
On quasi-uniform meshes, choosing $\delta\simeq h$ (equivalently $n\simeq h^{-d}\simeq\#\Th$) balances the three terms in the parenthesis. For $r\ge2$ this gives $\mathcal A(u)\lesssim h^{r-1}\|u\|_r$, and Theorem~\ref{thm:main} consequently guarantees an $O(h)$ error for sufficiently smooth solutions with bounded $|u|_2+\mu^{-1}|p|_1$. This is an upper-bound guarantee, not a lower bound on the achievable numerical error; higher observed convergence orders are not excluded.
\end{corollary}

\section{Numerical experiments}\label{sec:numerics}
In this section, we present numerical experiments that confirm our theoretical development. 

Unless stated otherwise, we take $\Omega=(-1,1)^2$, $\mu=\alpha=1$, $c_\tau=0.01$, ReLU$^3$ features with Fibonacci parameters on $\mathbb S^2$, and $N\times N$ squares, each cut into two triangles, so that $M=\#\Th=2N^2$ and the leg length is $h=2/N$. The network has no
boundary factor; the Dirichlet condition is imposed by the Nitsche terms of $A_h$. Here $n$ is the number of Fibonacci directions; neurons
whose hyperplane misses $\Omega$ are polynomials on $\Omega$ and are discarded, and the DoF reported in the tables count the network unknowns
after this step and after orthonormalization in $(\bseps\cdot,\bseps\cdot)$ plus the boundary penalty (relative cut $10^{-11}$), plus the $\#\Th$
enrichment unknowns. We also use $\tanh$ features $\phi_j=\tanh(a\,\theta_j\cdot(x,1))$ with $a=N/4$; they satisfy \eqref{eq:W}, so
Theorem~\ref{thm:main} applies, but Corollary~\ref{cor:rates} uses the approximation theory of ReLU$^k$ networks, so for $\tanh$ the rates are
observed, not proved. The exact solution is
\begin{eqnarray*}
u &=& \operatorname{curl}\Psi+\tfrac{\mathrm{Ra}}\lambda\nabla\phi,\quad p=\mathrm{Ra}\,\Delta\phi,\\
\Psi &=& \sin^2\pi x\,\sin^2\pi y,\quad \phi=(1+\cos\pi x)(1+\cos\pi y),
\end{eqnarray*}
so that $\mathrm{Ra}$ scales the gradient part of the load, and $u=0$ on $\partial\Omega$. Network integrals use a $7$-point rule on $16$
subtriangles per element and $3$-point Gauss rules on $4$ subsegments per face. The linear system is solved as the mixed system
\eqref{eq:method} by a direct solver, which also covers $\lambda=\infty$; for $\lambda<\infty$ it gives the solution of \eqref{eq:primal}. Errors
are $|u-U|_{1,h}$ and $\|p-P\|_0$ ($|u|_1=27.9$), relative where stated.

\begin{figure}[!htb]
\centering
\includegraphics[width=\textwidth]{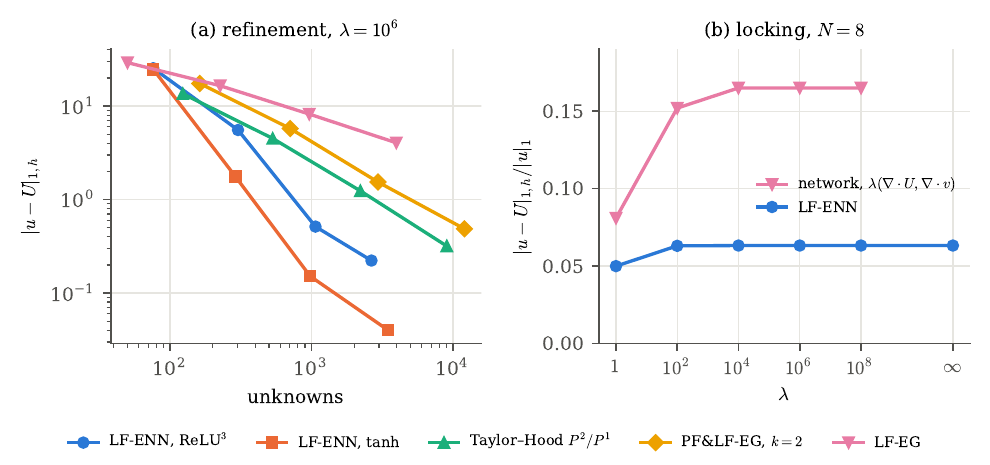}
\caption{(a) Error $|u-U|_{1,h}$ versus unknowns under refinement ($\lambda=10^6$, $n=\#\Th$; Tables~\ref{tab:order} and~\ref{tab:eg}).
(b) Relative error versus $\lambda$ for the network with the exact divergence and for LF-ENN ($N=8$, $n=256$; Table~\ref{tab:lambda}).}
\label{fig:conv2d}
\end{figure}

Table~\ref{tab:lambda} and Figure~\ref{fig:conv2d}(b) compare LF-ENN with the same network without the enrichment and with the exact divergence,
$\lambda(\nabla\cdot U,\nabla\cdot v)$ (with the same Nitsche and least-squares terms). The plain network loses accuracy as $\lambda$ grows, and its
pressure error grows to $222\%$, while the LF-ENN errors change only slightly for $\lambda\ge10^2$, including $\lambda=\infty$.

\begin{table}[!htb]
\centering
\caption{Relative error $|u-U|_{1,h}/|u|_1$ (relative pressure error in parentheses) versus $\lambda$; $N=8$, $n=256$, $\mathrm{Ra}=1$.
For $\lambda=\infty$ the plain network would have to be restricted to its exactly divergence-free subspace, which cannot be computed in a stable
way: with eigenvalue cut-offs $10^{-8}$, $10^{-10}$ and $10^{-12}$ of the divergence Gram matrix the error is $0.166$, $0.208$ and $0.402$.}
\label{tab:lambda}
\begin{tabular}{r|c|c}
\toprule
$\lambda$ & plain network, $\lambda(\nabla\cdot U,\nabla\cdot v)$ & LF-ENN\\
\midrule
$1$      & $0.0809$ ($0.233$) & $0.0499$ ($0.227$)\\
$10^2$   & $0.152$ ($0.848$)  & $0.0630$ ($0.259$)\\
$10^4$   & $0.165$ ($0.900$)  & $0.0632$ ($0.261$)\\
$10^6$   & $0.165$ ($0.991$)  & $0.0632$ ($0.261$)\\
$10^8$   & $0.165$ ($2.22$)   & $0.0632$ ($0.261$)\\
$\infty$ & --                 & $0.0632$ ($0.261$)\\
\bottomrule
\end{tabular}
\end{table}

Table~\ref{tab:order} and Figure~\ref{fig:conv2d}(a) refine the mesh with $n=\#\Th$. The observed orders are above
the guaranteed order one. With ReLU$^3$ features LF-ENN is $2.4$ and $1.4$ times more accurate than Taylor-Hood on the same mesh for
$N=16$ and $32$, with $48\%$ and $29\%$ of its unknowns; with $\tanh$ features it is $2.6$, $8.2$ and $7.9$ times more accurate for $N=8$, $16$ and $32$, with $55\%$, $44\%$ and $38\%$ of its unknowns.
On the coarsest meshes Taylor--Hood is more accurate. Along $n=\#\Th$ the computed inf-sup constants are $\beta=0.347$, $0.348$, $0.303$
(ReLU$^3$) and $0.391$, $0.320$, $0.260$ ($\tanh$) for $N=4$, $8$, $16$, consistent with Lemma~\ref{lem:infsup}; the relative pressure errors
are $0.84$, $0.39$, $0.12$ and $0.070$ (ReLU$^3$, $N=4,\dots,32$).

\begin{table}[!htb]
\centering
\caption{Refinement with $n=\#\Th$: $|u-U|_{1,h}$, observed orders in parentheses, and unknowns; $\lambda=10^6$, $\mathrm{Ra}=1$.
Taylor--Hood $P^2/P^1$ with the standard load; its DoF include the pressure.}
\label{tab:order}
\begin{tabular}{r|rl|rl|rl}
\toprule
 & \multicolumn{2}{c|}{LF-ENN, ReLU$^3$} & \multicolumn{2}{c|}{LF-ENN, $\tanh$} & \multicolumn{2}{c}{Taylor--Hood $P^2/P^1$}\\
$N$ & DoF & error & DoF & error & DoF & error\\
\midrule
4  & 76   & 25.42         & 76   & 24.76         & 123  & 13.57\\
8  & 302  & 5.550 (2.20)  & 290  & 1.757 (3.82)  & 531  & 4.509 (1.59)\\
16 & 1065 & 0.5138 (3.43) & 978  & 0.1512 (3.54) & 2211 & 1.239 (1.86)\\
32 & 2649 & 0.2225 (1.21) & 3459 & 0.0402 (1.91) & 9027 & 0.318 (1.96)\\
\bottomrule
\end{tabular}
\end{table}

The method has a locking-free error bound but is not pressure-robust (Theorem~\ref{thm:main}). For $\mathrm{Ra}=10^4$ ($\lambda=10^6$,
$n=\#\Th$) the ReLU$^3$ errors are $28\,583$, $6557$ and $1637$ for $N=4$, $8$, $16$ ($\tanh$: $29\,258$, $6262$, $1503$), against $22\,717$, $3319$
and $428$ for Taylor--Hood: both grow linearly in $\mathrm{Ra}$, and Taylor--Hood is less affected. Large gradient loads require the
pressure-robust variant of Remark~\ref{rem:choices}(c).

Table~\ref{tab:eg} compares with LF-EG \cite{YLZ22} (NIPG, penalty
$\beta_{\mathrm{EG}}=10^{-3}$, standard load) and with the parameter-free and locking-free EG method (PF\&LF-EG) of Su, Tong, Zhang and Zhang
\cite{STZZ24} with $k=1$ and $k=2$, in our implementations of their methods on the same meshes. For displacements vanishing on $\partial\Omega$
the weak divergence of PF\&LF-EG is our $\divh$ ($k=1$); the methods differ in the shear term, the enrichment for $k=2$, and the conforming space.
For $N\ge16$ LF-ENN (ReLU$^3$) is $3.0$ and $2.2$ times more accurate than PF\&LF-EG with $k=2$, with $36\%$ and $22\%$ of its unknowns, and with
$\tanh$ features $10$ and $12$ times for $N=16$ and $32$. All methods use the standard load, and for $\mathrm{Ra}=10^4$ the EG methods are less affected
(PF\&LF-EG, $k=2$: $14\,116$, $2131$, $497$ for $N=4$, $8$, $16$).

\begin{table}[!htb]
\centering
\caption{Comparison with enriched Galerkin methods on the same meshes; $\lambda=10^6$, $\mathrm{Ra}=1$: $|u-U|_{1,h}$ and unknowns (boundary
values eliminated). LF-ENN with $n=\#\Th$.}
\label{tab:eg}
\resizebox{\textwidth}{!}{%
\begin{tabular}{r|rc|rc|rc|rc|rc}
\toprule
 & \multicolumn{2}{c|}{LF-EG \cite{YLZ22}} & \multicolumn{2}{c|}{PF\&LF-EG, $k=1$} & \multicolumn{2}{c|}{PF\&LF-EG, $k=2$} & \multicolumn{2}{c|}{LF-ENN, ReLU$^3$} & \multicolumn{2}{c}{LF-ENN, $\tanh$}\\
$N$ & DoF & error & DoF & error & DoF & error & DoF & error & DoF & error\\
\midrule
4  & 50   & 29.20 & 50   & 28.07 & 162     & 17.49 & 76   & 25.42  & 76  & 24.76\\
8  & 226  & 16.52 & 226  & 15.56 & 706     & 5.749 & 302  & 5.550  & 290 & 1.757\\
16 & 962  & 8.22  & 962  & 8.006 & 2946    & 1.540 & 1065 & 0.5138 & 978 & 0.1512\\
32 & 3970 & 4.05  & 3970 & 4.031 & 12\,034 & 0.486 & 2649 & 0.2225 & 3459 & 0.0402\\
\bottomrule
\end{tabular}}
\end{table}

We clamp the side $x=-1$ and prescribe the exact traction $g=(2\mu\bseps(u)+p\mathbf I)\nu$ on the other
sides ($\lambda=10^6$, $n=\#\Th$, $c_H=0.01$). The errors are $25.15$, $5.683$ and $0.4984$ for $N=4$, $8$, $16$, the same to six digits for
$\lambda=\infty$ and within $1\%$ for $\lambda=10^2$ ($25.11$, $5.657$, $0.4953$); for $\mathrm{Ra}=10^2$ they are $292.9$, $64.85$ and $17.53$. Without the
Hessian term ($c_H=0$) the errors are $25.10$, $5.629$ and $0.4827$, so the term costs $0.2\%$, $1.0\%$ and $3.2\%$; with $c_H=0.1$ and $1$ the cost is
larger ($1.071$ and $3.616$ at $N=16$), as expected from its consistency error. Without the Hessian term the analysis would need the trace
constant of Remark~\ref{rem:CN}; it is $C_N=3.96$, $4.90$, $5.72$ for $N=4$, $8$, $16$ along $n=\#\Th$ and grows with $n$ on a fixed mesh
($8.46$ for $N=4$, $n=8\#\Th$).

\begin{figure}[!htb]
\centering
\includegraphics[width=\textwidth]{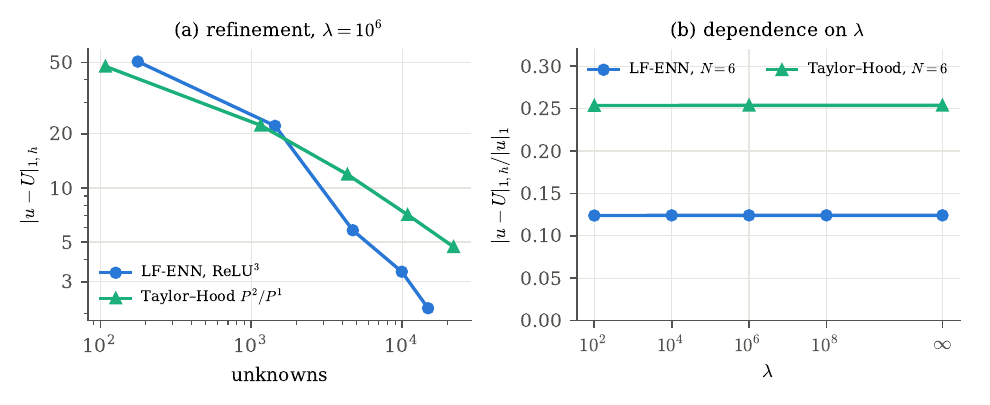}
\caption{Three dimensions. (a) $|u-U|_{1,h}$ versus unknowns ($\lambda=10^6$; Table~\ref{tab:3d}). (b) Relative error versus $\lambda$ for $N=6$.}
\label{fig:3d}
\end{figure}

On $\Omega=(-1,1)^3$ we divide $N^3$ cubes into six tetrahedra each (Kuhn's subdivision, $\#\Th=6N^3$), clamp
$\partial\Omega$, and take $u=\nabla\psi\times(1,1,1)^\top+\frac{\mathrm{Ra}}\lambda\nabla\varphi$, $p=\mathrm{Ra}\,\Delta\varphi$, with
$\psi=\prod_i\sin^2(\pi x_i)$ and $\varphi=\prod_i(1+\cos\pi x_i)$, the analogue of the planar solution; $\mathrm{Ra}=1$. LF-ENN uses ReLU$^3$
features $\sigma_3(\theta_j\cdot(x,1))$ with $n=\#\Th$ super-Fibonacci points $\theta_j$ on $\mathbb S^3$ \cite{Ale22}, the same parameters as in two
dimensions ($\alpha=1$, $c_\tau=0.01$, $h_e$ the face diameter, $h_T$ the element diameter), and no other change. Volume integrals use a
$27$-point collapsed Gauss rule on eight subtetrahedra, face integrals the $7$-point rule on four subtriangles. Table~\ref{tab:3d} and
Figure~\ref{fig:3d} compare with Taylor--Hood $P^2/P^1$ on the same meshes. On the two coarsest meshes the methods are comparable; for $N=8$
LF-ENN is $2.1$ times more accurate than Taylor--Hood with $92\%$ of its unknowns, and for $N=10$ it is $2.2$ times more accurate with $68\%$ of
its unknowns. For $N=10$ the network is orthonormalized componentwise, in the $H^1$-seminorm plus the boundary penalty, so that the dense
network matrices need not be stored; on the mesh $N=4$ this changes the error only in the sixth digit. The pressure of Taylor--Hood, which is continuous and piecewise linear, is more accurate. The LF-ENN
errors change by less than $0.3\%$ between $\lambda=10^2$ and $\lambda=\infty$ (for $N=6$: $5.800$, $5.810$, $5.810$, $5.810$, $5.810$ for
$\lambda=10^2,10^4,10^6,10^8,\infty$; for $N=10$: $2.133$, $2.138$, $2.138$ for $\lambda=10^2,10^6,\infty$), supporting locking-free behavior in the tested three-dimensional configurations, consistent with Theorem~\ref{thm:main}.

\begin{table}[!htb]
\centering
\caption{Three dimensions, clamped cube, $\lambda=10^6$, $n=\#\Th=6N^3$: $|u-U|_{1,h}$ (observed orders in parentheses), relative pressure error,
and unknowns. Taylor--Hood $P^2/P^1$ on the same meshes; its unknowns include the pressure.}
\label{tab:3d}
\begin{tabular}{r|rll|rll}
\toprule
 & \multicolumn{3}{c|}{LF-ENN, ReLU$^3$} & \multicolumn{3}{c}{Taylor--Hood $P^2/P^1$}\\
$N$ & DoF & error & pressure & DoF & error & pressure\\
\midrule
2  & 177  & 50.50        & 1.03  & 108     & 47.50        & 0.596\\
4  & 1437 & 22.16 (1.19) & 0.690 & 1154    & 22.31 (1.09) & 0.323\\
6  & 4706 & 5.810 (3.30) & 0.383 & 4336    & 11.88 (1.55) & 0.134\\
8  & 9936 & 3.415 (1.85) & 0.262 & 10\,854 & 7.090 (1.79) & 0.0655\\
10 & 14\,841 & 2.138 (2.10) & 0.201 & 21\,908 & 4.697 (1.85) & 0.0388\\
\bottomrule
\end{tabular}
\end{table}

\section{Cook's membrane}\label{sec:cook}
We use the setting of \cite{SGH19}: the panel with vertices $(0,0)$, $(48,44)$, $(48,60)$, $(0,44)$,
clamped at $x=0$, shear traction $g=(0,1/16)$ at $x=48$, $f=0$, $\mu=0.375$ and $\lambda=7.5\times10^6$ (nearly
incompressible). The quantity of interest is $u_y$ at $Q=(48,52)$, with reference $16.442$ \cite{SGH19}. LF-ENN uses the Nitsche terms on the
clamped side, $c_\tau=0.01$, the Hessian term with $c_H=0.01$ on the elements along the traction sides, $n=\#\Th$, and the network input mapped affinely to $[-1,1]^2$ (Section~\ref{sec:method}); the load is
$\langle g,v\rangle_{\Gamma_N}$. We compare with continuous $P^1$ (locks), $P^1+\Enr$ (Remark~\ref{rem:choices}(b)), and Taylor--Hood $P^2/P^1$
\cite{TH73,BF91}, on the same meshes. Unknowns are counted after orthonormalization and elimination of clamped values, including the pressure
for Taylor--Hood. The solution has a corner singularity where the clamped side meets the upper boundary. As discussed in Remark~\ref{rem:corner}, applying Theorem~\ref{thm:lowreg} requires verification of the relevant corner regularity exponent.

\begin{figure}[!htb]
\centering
\includegraphics[width=\textwidth]{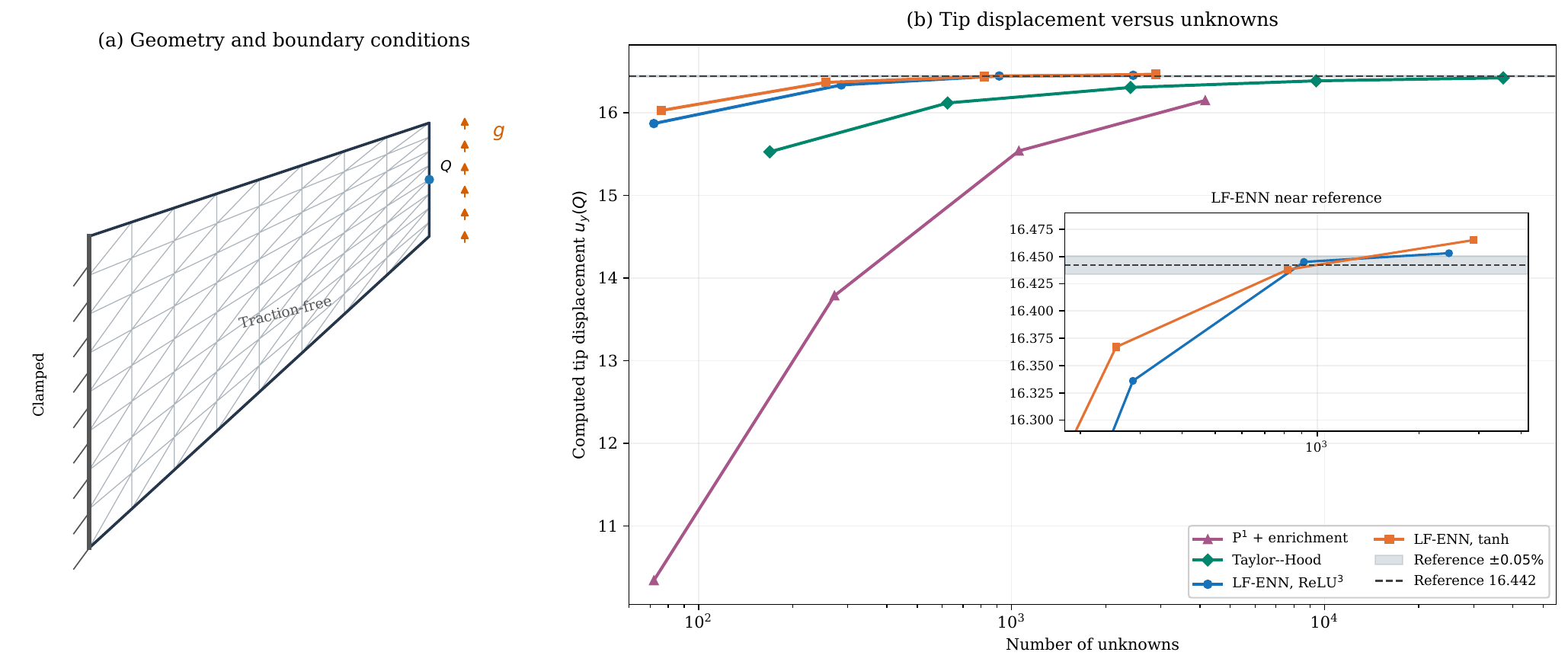}
\caption{Cook's membrane. (a) Geometry, boundary conditions, and the structured mesh with $N=8$. (b) Computed vertical tip displacement $u_y(Q)$ versus the number of unknowns for the methods in Table~\ref{tab:cook}. The dashed horizontal line indicates the reported reference value $u_y(Q)=16.442$; the light shaded band denotes the reported approximate $\pm0.05\%$ reference uncertainty (not a rigorous confidence interval). The inset enlarges the LF-ENN results near this reference. The locking $P^1$ method is omitted from panel (b).}
\label{fig:cook}
\end{figure}

Table~\ref{tab:cook} reports the computed tip displacement $u_y(Q)$, the corresponding relative errors, and the numbers of unknowns, while Figure~\ref{fig:cook}(b) plots $u_y(Q)$ against the number of unknowns. The continuous $P^1$ method remains at approximately $28\%$ of the reference displacement, and $P^1+\Enr$ converges slowly. Taylor--Hood requires $37\,249$ unknowns to attain a relative error of $0.12\%$, whereas LF-ENN agrees with the reported reference value within $0.02\%$ using $914$ unknowns (ReLU$^3$) and within $0.025\%$ using $819$ unknowns ($\tanh$).

The LF-ENN values are independent of $\lambda$ to five digits for $\lambda$ from $7.5\times10^4$ to $7.5\times10^{10}$ and for $\lambda=\infty$ (for example, $16.3356$ at $N=8$). Without the Hessian term, the values differ by at most $0.6\%$ ($15.951$, $16.350$, $16.447$, and $16.454$ for ReLU$^3$).

As shown in Figure~\ref{fig:cook}(b), the LF-ENN tip displacements increase under refinement. The ReLU$^3$ results cross the reported reference value $16.442$ between $N=8$ and $N=16$, whereas the $\tanh$ results cross it between $N=16$ and $N=32$. Consequently, the absolute relative errors are nonmonotone and, for both networks, are smaller at $N=16$ than at $N=32$. This does not necessarily indicate deterioration of the global approximation, since the comparison concerns a single pointwise displacement rather than an energy-norm error. Moreover, the reported reference value has an uncertainty of approximately $0.05\%$, shown as a shaded band in Figure~\ref{fig:cook}(b), comparable to the smallest observed relative errors. This band represents a reported tolerance rather than a statistical confidence interval.
\begin{table}[!htb]
\centering
\caption{Cook's membrane, structured meshes: $u_y(Q)$ (relative error) and unknowns; reference $16.442$.}
\label{tab:cook}
\resizebox{\textwidth}{!}{%
\begin{tabular}{r|cr|cr|cr|cr|cr}
\toprule
 & \multicolumn{2}{c|}{$P^1$} & \multicolumn{2}{c|}{$P^1+\Enr$} & \multicolumn{2}{c|}{Taylor--Hood $P^2/P^1$} & \multicolumn{2}{c|}{LF-ENN, ReLU$^3$} & \multicolumn{2}{c}{LF-ENN, $\tanh$}\\
$N$ & $u_y(Q)$ & DoF & $u_y(Q)$ & DoF & $u_y(Q)$ & DoF & $u_y(Q)$ & DoF & $u_y(Q)$ & DoF\\
\midrule
4  & 4.605 (72.0\%) & 40   & 10.341 (37.1\%) & 72   & 15.527 (5.57\%) & 169    & 15.869 (3.49\%) & 72   & 16.027 (2.52\%) & 76\\
8  & 4.630 (71.8\%) & 144  & 13.786 (16.2\%) & 272  & 16.118 (1.97\%) & 625    & 16.336 (0.65\%) & 286  & 16.367 (0.46\%) & 255\\
16 & 4.636 (71.8\%) & 544  & 15.538 (5.50\%) & 1056 & 16.306 (0.83\%) & 2401   & 16.445 (0.02\%) & 914  & 16.438 (0.025\%) & 819\\
32 & 4.638 (71.8\%) & 2112 & 16.149 (1.78\%) & 4160 & 16.386 (0.34\%) & 9409   & 16.453 (0.07\%) & 2449 & 16.465 (0.14\%) & 2894\\
64 & --             & --   & --              & --   & 16.422 (0.12\%) & 37\,249 & -- & -- & -- & --\\
\bottomrule
\end{tabular}}
\end{table}

On smoothed random Delaunay meshes with $92$, $344$ and $1323$ triangles (three meshes each, minimum angles
$15^\circ$--$26^\circ$) LF-ENN (ReLU$^3$) gives $1.2$--$1.9\%$, $0.61$--$0.68\%$ and $0.12$--$0.18\%$ with $208$, $695$ and $1768$ unknowns. On random
Delaunay meshes of comparable size ($91$--$93$, $365$--$377$, $1477$--$1485$ triangles) Taylor--Hood gives $0.51$--$0.55\%$, $0.17$--$0.22\%$ and
$0.02$--$0.25\%$ with about $460$, $1760$ and $6850$ unknowns. On unstructured meshes Taylor--Hood is thus more accurate on the coarse meshes, and
the two methods are comparable on the finest ones, LF-ENN with a quarter of the unknowns. The FCFV method of \cite{SGH19} is reported within
$0.5\%$ with $2.6\times10^5$ equations.

\section{Discussion and open problems}\label{sec:conclusions}
LF-ENN is a displacement method in which a linearized neural network becomes locking-free with one enrichment unknown per element, because the enrichment completes every conforming field to a discretely divergence-free one (Proposition~\ref{prop:kernel}). The constraint is local
conservation of volume, $\int_{\partial T\setminus\Gamma_D}\avg U\cdot\nu=|T|P_T/\lambda$, with the element pressure as multiplier. For each fixed shape-regular mesh the error estimate is independent of $\lambda$, although its constants may depend on the discrete inf-sup constant. The proven upper bound includes an $O(h)$ consistency term, so increasing network width alone does not improve the guaranteed order; this does not rule out smaller errors in practice. This motivates the pairing $n\approx\#\Th$ on quasi-uniform meshes. The Nitsche terms,
the least-squares term and the Hessian term along the traction boundary remove the boundary factor and every trace constant of the
network, so that, when the network resolves the mesh, all constants are independent of $n$, $h$ and $\lambda$. The price is that the method is not pressure-robust, and
the reported implementation uses dense network blocks, which can be more costly than sparse Taylor--Hood systems; a controlled runtime and memory comparison remains desirable.

There are several open issues and work in progress. First, a Bernstein inequality for linearized networks is still an open question, with which the Hessian term will be unnecessary. Secondly, an explanation of the observed orders above one will need to be established. Thirdly, variants with discontinuous $P^m$ pressure could be designed, in which the least-squares term can
include the pressure gradient and the consistency order is raised. Lastly, $\lambda$-robust preconditioning and a sparse implementation could be developed. 

\section*{CRediT authorship contribution statement}
\textbf{Bonyong Koo:} Conceptualization, Software, Validation, Investigation, Writing- review \& editing, Funding acquisition.
\textbf{Young Ju Lee:} Conceptualization, Methodology, Formal analysis, Software, Writing -original draft,
Writing -- review \& editing, Supervision. 

\section*{Declaration of competing interest}
The authors declare that they have no known competing financial interests or personal relationships that
could have appeared to influence the work reported in this paper.

\section*{Data availability}
The source code and input data used for the numerical experiments are not included in this submission. The authors should provide an availability statement or repository information before publication.

\section*{Acknowledgments}
Y.-J.~Lee was supported in part by the U.S. National Science Foundation under award DMS-2208499. B. Koo was supported by the National Research Foundation of Korea (NRF) grant by the Korean government (No. RS-2026-25598463).

\bibliographystyle{elsarticle-num}
\bibliography{enn_lf}
\end{document}